\documentclass[11pt]{article}
\usepackage{amssymb,amsmath,mathrsfs,txfonts,graphicx,color,tikz}
\graphicspath{{figs/}}

\begin{document}

\allowdisplaybreaks
\def\pd#1#2{\frac{\partial#1}{\partial#2}}
\def\dfrac{\displaystyle\frac}
\let\oldsection\section
\renewcommand\section{\setcounter{equation}{0}\oldsection}
\renewcommand\thesection{\arabic{section}}
\renewcommand\theequation{\thesection.\arabic{equation}}

\newtheorem{theorem}{\indent Theorem}[section]
\newtheorem{lemma}{\indent Lemma}[section]
\newtheorem{proposition}{\indent Proposition}[section]
\newtheorem{definition}{\indent Definition}[section]
\newtheorem{remark}{\indent Remark}[section]
\newtheorem{corollary}{\indent Corollary}[section]

\title{\LARGE
Sharp, Smooth, and Oscillatory Traveling waves of Degenerate Diffusion Equation with Delay}
\author{
Tianyuan Xu$^{a,d}$, Shanming Ji$^{b,d}$\thanks{Corresponding author, email:jism@scut.edu.cn},
Ming Mei$^{c,d}$, Jingxue Yin$^a$,
\\
\\
{ \small \it $^a$School of Mathematical Sciences, South China Normal University}
\\
{ \small \it Guangzhou, Guangdong, 510631, P.~R.~China}
\\
{ \small \it $^b$School of Mathematics, South China University of Technology}
\\
{ \small \it Guangzhou, Guangdong, 510641, P.~R.~China}
\\
{ \small \it $^c$Department of Mathematics, Champlain College Saint-Lambert}
\\
{ \small \it Quebec,  J4P 3P2, Canada, and}
\\
{ \small \it $^d$Department of Mathematics and Statistics, McGill University}
\\
{ \small \it Montreal, Quebec,   H3A 2K6, Canada}
}
\date{}

\maketitle

\begin{abstract}
We consider the non-monotone degenerate diffusion equation with time delay.
Different from the linear diffusion equation, the degenerate equation allows for
semi-compactly supported traveling waves.
In particular, we discover sharp-oscillating waves with sharp edges and non-decaying oscillations.
The degenerate diffusion and the effect of time delay cause us essential difficulties.
We show the existence for both sharp and smooth traveling wave solutions.
Furthermore, we prove the oscillating properties of the waves for large
wave speeds and large time delay.
Since the existing approaches are not applicable, we develop
a new technique to show the existence of the sharp, smooth and oscillatory traveling waves.
\end{abstract}


\section{Introduction}
In this paper, we are concerned with the
traveling wave solutions to a degenerate diffusion equation
with time delay
\begin{align}\label{eq-main}
\begin{cases}
\displaystyle
\pd u t=D\Delta u^m-d(u)+
b(u(t-r,x)),\quad &x\in \mathbb R,~ t>0,\\
u(s,x)=u_0(s,x), &x\in \mathbb R,~ s\in[-r,0],
\end{cases}
\end{align}
which models the population dynamics for single species with age structure.
Here,
$D$ denotes the diffusion coefficient,
$u$ represents the density of total mature population
at location $x$ and time $t$,
$D\Delta u^m$ is
the density-dependent diffusion.
Such a degenerate diffusion means that the spatial-diffusion rate increases with  population density, particularly, zero
density implies non-diffusion.
This is with more ecological sense \cite{Carl,Murry,Okubo}. Two nonlinear functions
$b(u)$ and $d(u)$ represent the birth rate
and the death rate of the matured respectively.

From biological experiments,
\eqref{eq-main} admits two constant equilibria $u_-=0$ and $u_+=\kappa>0$,
where $u_-=0$ is unstable and
$u_+=\kappa$ is stable for the
spatially homogeneous equation associated with \eqref{eq-main}.
Our model includes the classical Fisher-KPP
equation \cite{Fisher,Gourley00} and a lots of evolution equations in ecology,
for example,
the well-studied diffusive Nichloson's blowflies equation and Mackey-Glass  equation \cite{Faria06,LiJNS,LinMei,Mei_LinJDE09,Mei04} with
the death function
$d(u)=\delta u$, the birth function
$$b_1(u)=pue^{-au^q}, \mbox{ or }  b_2(u)=\frac{pu}{1+au^q},
\quad p>0,~ q>0,~ a>0;$$
and the age-structured population model \cite{Chern,Gourley04,Kuang03,Li_Mei} with
$$d(u)=\delta u^2, \mbox{ and } b(u)=pe^{-\gamma r}u,
\quad p>0,~ \delta>0,~ \gamma>0.$$

Our main purpose is to study the existence and non-existence of both sharp and smooth traveling wave solutions together with the
oscillatory properties for the system \eqref{eq-main} without the
monotonicity assumption on the birth function $b(\cdot)$.
A traveling wave solution is
a specific form of solutions with $u(t,x)=\phi(x+ct)$, where $c$ is the wave speed.
The ecological meaning of traveling wave solutions for \eqref{eq-main} is that
the individuals of the population disperse throughout the habitat in a certain density profile
moving with a constant speed.

The study of the invasion and spreading of species with linear diffusion has
a long history.
Since the pioneering work of Schaaf \cite{Schaaf}, the existence of traveling waves to
reaction diffusion equations has been extensively studied.
The authors
So, Wu and Zou \cite{Joseph} proved the existence of monotone traveling wave solutions by the upper and lower solutions method.
Faria and Trofumichunk \cite{Faria06} found that the traveling waves can be oscillatory
when the time delay is large.
Using fixed points method, Ma \cite{MaJDE} proved the existence of non-monotone traveling waves for the time-delayed equation with nonlocal birth rate function.
Gomez and Trofimchunk \cite{Gomez} proved the existence of oscillatory
and monotone traveling waves for any time delay.
Alfaro et al. \cite{Alfaro} combined the priori estimates and the
Leray Schauder topology method to study the existence of oscillatory traveling waves
for the non-monotone bistable equation with time delay.
In \cite{Oujde10}, a new approach based on the shooting method with upper and lower solutions
was developed to study the time-delayed Fisher-KPP
equation with non-monotone source.
The global stability of critical traveling waves
with optimal decay estimates are investigated in \cite{IJNAM}.

The first application of reaction-diffusion equation in biology was to use linear diffusion to model spatial diffusion of population \cite{Fisher,Skellam}.
There are, however, considerable evidences that
several species migrate from densely populated areas to sparsely areas to avoid overcrowding, rather than random walk diffusion \cite{Carl,Morisita1971}.
Gurney and Nisbet \cite{Gurney75} first proposed density-dependent dispersal to describe population spreading.
This positive density-dependent mechanism arises from
competition between conspecifics or deteriorating environmental conditions \cite{Matthysen}.
Now it is a common feature of population spreading modelling in ecology.

Dynamical behaviors of traveling wave for  degenerate reaction-diffusion equation are extremely
rich and interesting.
The degeneracy raises the possibility of sharp type traveling waves.
Different from the smooth traveling waves, in the sharp type waves, the population density $u$ decreases to zero at a finite point, rather than decaying to zero asymptotically.
Sharp traveling waves are sometimes called finite waves.
In 1980s,  Aronson \cite{Aronson80} first studied the sharp waves with critical wave speed for degenerate diffusion equation without time delay.
Then Pablo and Vazqueze \cite{Pablo} found the sharp waves for more degenerate Fisher-KPP equations.
In 2005, Gilding and Kersner \cite{Gilding} obtained the exact sharp waves for a particular
Fisher-KPP equation with degenerate diffusion and convection.
Recently, some detailed discussions of degenerate diffusion with time delay are emerging.
Huang et al. \cite{Huang-Jin-Mei-Yin} first obtained the existence and stability of time-delayed population dynamics model with degenerate diffusion for small time delay.
Later then,
we \cite{JDE18XU} proved the existence of monotone traveling wave solutions for large time delay.
The approach adapted for the
proof is the monotone technique with the viscosity vanishing method.

The two most important questions in biological spreading processes ask how
fast the population propagates and what shape it forms.
In this paper, we work on the traveling waves for a non-monotone degenerate reaction
diffusion equation with large time delay.
We focus on the influence of the diffusion and the non-monotone birth rate function on the existence and shape of such profiles.
The wave behavior is rather complicated and rich for the degenerate diffusion
equations with time delay.
There exist smooth traveling waves,
sharp waves, and both of these waves show big oscillations for large
wave speeds and large time delay.

We first prove the existence of smooth traveling wave solutions for model \eqref{eq-main}.
The wave profile equations are usually solved either
through the iteration procedure or by means of the phase plane analysis.
These approaches lead to restrictive assumptions such as monotonicity or small time delay on the delayed term.
Our problem does not admit any comparison principle but possesses large time delay.
This prevents the application of classical techniques, and we need to introduce new ideas and techniques to overcome the emerging difficulties caused by large time delay and non-monotonicity
as well as the degeneracy of diffusion.
Using the Schauder Fixed Point Theorem,  we construct
an appropriate profile set  with upper and lower profiles
for two auxiliary problems and obtained the existence of  monotone and non-monotone traveling waves.

The time delay $r$ and  the degenerate diffusion  in model \eqref{eq-main} have a strong influence on the geometry of wave profiles,
such as sharp waves caused by  the degeneracy.
We emphasize that it is the first literature on the existence of sharp type traveling wave with time delay as far as we known.
A sharp wave solution $\phi(t)$ is a special solution with semi-compact support such that
$\phi(t)\equiv0$ for $t\le t_0$ and $\phi(t)>0$ for $t>t_0$ with some $t_0\in\mathbb R$.
The existence of sharp wave solutions for the case without time delay
and with Nicholson's birth rate function $b(u)=pue^{-au}$ and death rate function
$d(u)=\delta u$ for some constants $p,a,\delta$ is proved in \cite{Huang-Jin-Mei-Yin}.
Due to the lack of monotonicity  and the
bad effect of time delay, the method of \cite{Aronson80,Pablo,Gilding,Huang-Jin-Mei-Yin} are not applicable.
Based on an observation of the delicate structure of time delay and sharp edge,
a new delayed iteration approach is developed to solve the delayed degenerate equation.
To our best knowledge, this is the first framework
of showing the existence of sharp traveling wave solution for the degenerate
diffusion equation with large time delay.

The speed selection mechanism for the degenerate diffusion equation with time delay \eqref{eq-main} is interesting.
We prove the nonexistence of traveling waves for sub-critical wave speed.
The proof is based on the phase transform approach
with some modification suitable for large time delay
and non-monotone birth rate functions. The critical wave speeds for
both the sharp waves and smooth waves of model \eqref{eq-main} are nonlinearly determined.
The appearance of degenerate diffusion leads to the
failure of the ``linear determinacy principle'' \cite{Lewis}.
The wave behaviors cannot be determined by the linearization around
equilibrium zero but controlled by the whole wave structure.

Finally,  we investigate the oscillatory properties of the traveling waves in $+\infty$:
convergence to the positive equilibrium $\kappa$ when the delay and wave speed are small,
whereas oscillations around $\kappa$ for both large delay and large wave speed.
We give an explicit description of wave behaviors, depending on
the properties of the birth rate function, the tails of the waves may
approach the carrying capacity monotonically, may approach the carrying capacity $\kappa$ in an oscillatory manner, or may oscillate infinitely
around the carrying capacity, where its values are bounded above and below.

In the degenerate diffusion equation with time delay,
propagating traveling waves may possess different dynamical behaviors.
As far as we know, the study of wave profiles done in this paper is new and our results can be derived by none of the papers quoted above.
On the one hand,
Theorem \ref{th-semifinite}, Theorem \ref{th-osc} and Theorem \ref{th-nondecay}
presented later in this paper
imply that oscillating traveling waves with sharp type leading edge
(see Figure \ref{fig-sharp} and Figure \ref{fig-C1} for illustration) can appear.
Here, we call this type of traveling waves  ``sharp-oscillatory waves''.
On the other hand, we give a precise characterization of the geometric
dynamics of traveling waves.
Actually, as shown in Figure \ref{fig-cr} and Figure \ref{fig-cr-2},
the shapes of wave fronts can be predicted by the degenerate diffusion equation with time delay \eqref{eq-main}.
Patterns can be characterized by their velocity of propagation $c$, time delay $r$ and the degeneracy
index $m$.

The rest of this paper is organized as follows.
In Section 2, we present the main results on the existence and nonexistence
of traveling waves and the oscillatory properties of traveling waves.
Section 3 is devoted to the proof of the existence of the non-monotone smooth
traveling wave solutions, while in Section 4 we prove the nonexistence of traveling waves.
The existence of sharp traveling waves is proved in Section 5.
Finally, the oscillation properties of
traveling wave solutions are investigated in Section 6.

\section{Main results}
We consider the initial-value problem \eqref{eq-main},
where the time delay $r\ge0$, $m>1$, $D>0$,
$u_0\in L^2((-r,0)\times\Omega)$ for any compact set $\Omega\subset \mathbb R$.
Since \eqref{eq-main} is degenerate for $u=0$, we employ
the following definition of weak solutions.

\begin{definition} \label{de-weak}
A function $u\in L_\mathrm{loc}^2((0,+\infty)\times\mathbb R)$
is called a weak solution of \eqref{eq-main} if
$0\le u\in L^\infty((0,+\infty)\times\mathbb R)$,
$\nabla u^m\in L_{\mathrm{loc}}^2((0,+\infty)\times\mathbb R)$,
and for any $T>0$ and $\psi\in C_0^\infty((-r,T)\times\mathbb R)$
\begin{align*}
&-\int_0^T\int_{\mathbb R}u(t,x)\pd{\psi}{t}dxdt
+D\int_0^T\int_{\mathbb R}\nabla u^m\cdot\nabla\psi dxdt
+\int_0^T\int_{\mathbb R}d(u(t,x))\psi dxdt \\
&=\int_{\mathbb R}u_0(0,x)\psi(0,x)dx+
\int_r^{\max\{T,r\}}\int_{\mathbb R}
b(u(t-r,x))\psi(x,t)dxdt\\
&\ \ \ +\int_0^{\min\{T,r\}}\int_{\mathbb R}
b(u_0(t-r,x))\psi(x,t)dxdt.
\end{align*}
\end{definition}

We are looking for traveling wave solutions of \eqref{eq-main}
connecting the two equilibria $0$ and $\kappa>0$ in some sense that
they may oscillate around the positive equilibrium $\kappa$.
Let $\phi(\xi)$ with $\xi=x+ct$ and $c>0$
be a traveling wave solution of \eqref{eq-main},
we get (we write $\xi$ as $t$ for the sake of simplicity)
\begin{equation}\label{eq-tw}
c\phi'(t)=D({\phi^m}(t))''-d(\phi(t))
+b(\phi(t-cr)), \qquad t\in\mathbb R.
\end{equation}

The wave solution $\phi(t)$ may be non-monotone and even non-decaying oscillating around
the positive equilibrium $\kappa$ since the birth function $b(u)$ is non-monotone.
Meanwhile, it is also expected that the degenerate diffusion
equation \eqref{eq-main} may admit sharp type wave solution with semi-compact support.
So let us fix some terminology before going further.

\begin{definition} \label{de-tw}
A function $0\le\phi(t)\in W_{\mathrm{loc}}^{1,1}(\mathbb R)\cap L^\infty(\mathbb R)$
with $\phi^m(t)\in W_{\mathrm{loc}}^{1,1}(\mathbb R)$
is said to be a semi-wavefront of \eqref{eq-main} if
\\ \indent
(i) the profile function $\phi$ satisfies \eqref{eq-tw} in the sense of distributions,
\\ \indent
(ii) $\phi(-\infty)=0$, and
$0<\liminf_{t\to+\infty}\phi(t)\le\limsup_{t\to+\infty}\phi(t)<+\infty$,
\\ \indent
(iii) the leading edge of $\phi(t)$ near $-\infty$ is
monotonically increasing in the sense
that there exists a maximal interval $(-\infty,t_0)$ with $t_0\in(-\infty,+\infty]$
such that $\phi(t)$ is monotonically increasing in it
and if $t_0<+\infty$ then $\phi(t_0)>\kappa$.
We say that $t_0$ is the boundary of the leading edge of $\phi$.

A semi-wavefront $\phi(t)$ is said to be a wavefront of \eqref{eq-main}
if $\phi$ converges to $\kappa$ as $t$ tends to $+\infty$, i.e., $\phi(+\infty)=\kappa$.

A semi-wavefront (including wavefront) is said to be sharp if
there exists a $t_*\in\mathbb R$ such that $\phi(t)=0$ for all $t\le t_*$
and $\phi(t)>0$ for all $t>t_*$.
Otherwise, it is said to be a smooth semi-wavefront (or smooth wavefront)
if $\phi(t)>0$ for all $t\in\mathbb R$.

Furthermore, for the sharp semi-wavefronts (including wavefronts) $\phi(t)$,
if $\phi''\not\in L_\text{loc}^1(\mathbb R)$,
we say that $\phi(t)$ is a non-$C^1$ type sharp waves;
otherwise, if $\phi''\in L_\text{loc}^1(\mathbb R)$,
we say that $\phi(t)$ is a $C^1$ type sharp waves.
\end{definition}

\begin{figure}[htb]
\begin{center}
\begin{tikzpicture}[scale=0.85,domain=-7:7]
\def\axll{-7} \def\axlr{-3} \def\axml{-2} \def\axmr{2} \def\axrl{3} \def\axrr{7}
\def\ay{2} \def\ayt{3} \def\pa{\ay/180} \def\pb{8} \def\pc{18} \def\pd{8}
\draw[->,>=latex,line width=.5pt] (\axll,0)--(\axlr,0) node[right] {$\xi$};
\draw[dashed, line width=.5pt] (\axll,\ay)--(\axlr,\ay);
\draw[color=blue,line width=.7pt]
    plot[domain=-7:-3, samples=144, smooth] (\x,{\ay/2+\pa*atan(\pb*(\x+5))});
\node at (-5,-.4) {(A1)};
\draw[->,>=latex,line width=.5pt] (\axml,0)--(\axmr,0) node[right] {$\xi$};
\draw[dashed, line width=.5pt] (\axml,\ay)--(\axmr,\ay);
\draw[color=blue,line width=.7pt]
    plot[domain=-2:0, samples=144, smooth]
    (\x,{\ay/2+\pa*atan(\pb*(\x))});
\draw[color=blue,line width=.7pt]
    plot[domain=0:2, samples=144, smooth]
    (\x,{\ay/2+\pa*atan(\pb*(\x))+8*\x/(1+\pc*\x*\x+\pc*\x^4)*sin(\x*180*2+\x^2*180+\x^3*180)});
\node at (0,-.4) {(A2)};
\draw[->,>=latex,line width=.5pt] (\axrl,0)--(\axrr,0) node[right] {$\xi$};
\draw[dashed, line width=.5pt] (\axrl,\ay)--(\axrr,\ay);
\draw[color=blue,line width=.7pt]
    plot[domain=3:5, samples=144, smooth] (\x,{\ay/2+\pa*atan(\pb*(\x-5))});
\draw[color=blue,line width=.7pt]
    plot[domain=5:7, samples=144, smooth]
    (\x,{\ay/2+\pa*atan(\pb*(\x-5))+sin((\x-5)^2*180*2)});
\node at (5,-.4) {(A3)};
\end{tikzpicture}
\end{center}
\caption{Smooth type traveling waves: (A1) monotone wavefront; (A2) non-monotone wavefront;
(A3) non-decaying oscillating semi-wavefront.}
\label{fig-smooth}
\end{figure}
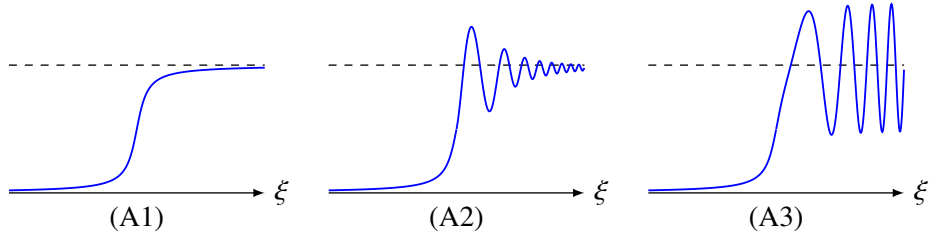

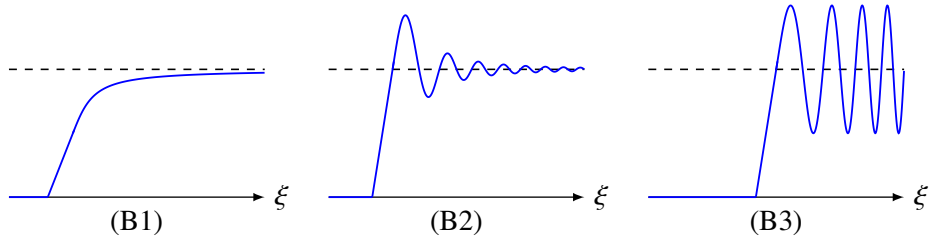
\begin{figure}[htb]
\begin{center}
\begin{tikzpicture}[scale=.85,domain=-7:7]
\def\axll{-7} \def\axlr{-3} \def\axml{-2} \def\axmr{2} \def\axrl{3} \def\axrr{7}
\def\ay{2} \def\ayt{3} \def\pa{\ay/180} \def\pb{4} \def\pc{12} \def\pd{8}
\def\ra{0.3927} \def\rb{8/3.14159} \def\sa{0.31831} \def\sb{2*3.14159}
\draw[->,>=latex,line width=.5pt] (\axll,0)--(\axlr,0) node[right] {$\xi$};
\draw[dashed, line width=.6pt] (\axll,\ay)--(\axlr,\ay);
\draw[color=blue,line width=.7pt] (\axll,0)--(-6.3927,0);
\draw[color=blue,line width=.7pt]
    plot[domain=-6.3927:-5.98, samples=144, smooth] (\x,{\rb*(\x+6.3927)});
\draw[color=blue,line width=.7pt]
    plot[domain=-6:-3, samples=144, smooth] (\x,{\ay/2+\pa*atan(\pb*(\x+6))});
\node at (-5,-.4) {(B1)};
\draw[->,>=latex,line width=.5pt] (\axml,0)--(\axmr,0) node[right] {$\xi$};
\draw[dashed, line width=.6pt] (\axml,\ay)--(\axmr,\ay);
\draw[color=blue,line width=.7pt] (\axml,0)--(-1.31831,0);
\draw[color=blue,line width=.7pt]
    plot[domain=-1.31831:-.98, samples=144, smooth] (\x,{\sb*(\x+1.31831)});
\draw[color=blue,line width=.7pt]
    plot[domain=-1:2, samples=144, smooth]
    (\x,{\ay+sin((\x+1)*180*2+(\x+1)^2*180)/(1+(\x+1)^2*4)});
\node at (0,-.4) {(B2)};
\draw[->,>=latex,line width=.5pt] (\axrl,0)--(\axrr,0) node[right] {$\xi$};
\draw[dashed, line width=.6pt] (\axrl,\ay)--(\axrr,\ay);
\draw[color=blue,line width=.7pt] (\axrl,0)--(4.68169,0);
\draw[color=blue,line width=.7pt]
    plot[domain=4.68169:5.02, samples=144, smooth] (\x,{\sb*(\x-4.68169)});
\draw[color=blue,line width=.7pt]
    plot[domain=5:7, samples=144, smooth]
    (\x,{\ay+sin((\x-5)*180*2+(\x-5)^2*180)});
\node at (5,-.4) {(B3)};
\end{tikzpicture}
\end{center}
\caption{Sharp type traveling waves --- non-$C^1$ type:
(B1) monotone wavefront; (B2) non-monotone wavefront;
(B3) non-decaying oscillating semi-wavefront.}
\label{fig-sharp}
\end{figure}

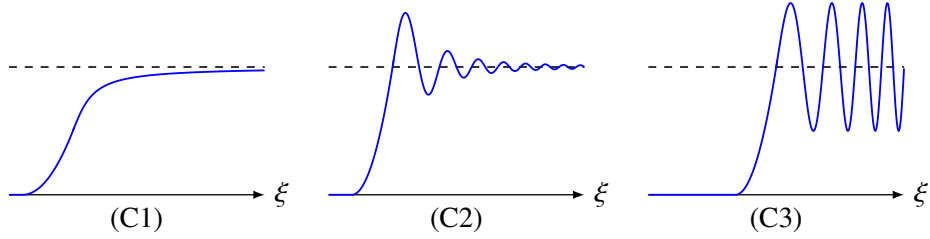
\begin{figure}[htp]

\begin{center}
\begin{tikzpicture}[scale=0.85,domain=-7:7]
\def\axll{-7} \def\axlr{-3} \def\axml{-2} \def\axmr{2} \def\axrl{3} \def\axrr{7}
\def\ay{2} \def\ayt{3} \def\pa{\ay/180} \def\pb{4} \def\pc{12} \def\pd{8}
\def\ra{0.3927} \def\rb{8/3.14159} \def\sa{0.31831} \def\sb{2*3.14159}
\draw[->,>=latex,line width=.5pt] (\axll,0)--(\axlr,0) node[right] {$\xi$};
\draw[dashed, line width=.6pt] (\axll,\ay)--(\axlr,\ay);
\draw[color=blue,line width=.7pt] (\axll,0)--(-6.7854,0);
\draw[color=blue,line width=.7pt]
    plot[domain=-6.7854:-5.98, samples=144, smooth]
    (\x,{(\x+6.7854)^2/.7854^2});
\draw[color=blue,line width=.7pt]
    plot[domain=-6:-3, samples=144, smooth]
    (\x,{\ay/2+\pa*atan(\pb*(\x+6))});
\node at (-5,-.4) {(C1)};
\draw[->,>=latex,line width=.5pt] (\axml,0)--(\axmr,0) node[right] {$\xi$};
\draw[dashed, line width=.6pt] (\axml,\ay)--(\axmr,\ay);
\draw[color=blue,line width=.7pt] (\axml,0)--(-1.63662,0);
\draw[color=blue,line width=.7pt]
    plot[domain=-1.63662:-.98, samples=144, smooth]
    (\x,{2*(\x+1.63662)^2/.63662^2});
\draw[color=blue,line width=.7pt]
    plot[domain=-1:2, samples=144, smooth]
    (\x,{\ay+sin((\x+1)*180*2+(\x+1)^2*180)/(1+(\x+1)^2*4)});
\node at (0,-.4) {(C2)};
\draw[->,>=latex,line width=.5pt] (\axrl,0)--(\axrr,0) node[right] {$\xi$};
\draw[dashed, line width=.6pt] (\axrl,\ay)--(\axrr,\ay);
\draw[color=blue,line width=.7pt] (\axrl,0)--(4.36338,0);
\draw[color=blue,line width=.7pt]
    plot[domain=4.36338:5.02, samples=144, smooth]
    (\x,{2*(\x-4.36338)^2/.63662^2});
\draw[color=blue,line width=.7pt]
    plot[domain=5:7, samples=144, smooth]
    (\x,{\ay+sin((\x-5)*180*2+(\x-5)^2*180)});
\node at (5,-.4) {(C3)};
\end{tikzpicture}
\end{center}
\caption{Sharp type traveling waves --- $C^1$ type:
(C1) monotone wavefront; (C2) non-monotone wavefront;
(C3) non-decaying oscillating semi-wavefront.}
\label{fig-C1}
\end{figure}

According to the above definition, the possible traveling wave solutions
are classified into monotone wavefronts, non-monotone wavefronts, or
non-decaying oscillating semi-wavefronts considering the monotonicity near $+\infty$;
and at the same time these waves can also be classified into sharp
or smooth type concerning the degeneracy near $-\infty$.
Moreover, the sharp type waves are further classified into $C^1$ type and non-$C^1$ type
according to the regularity.
See Figure \ref{fig-smooth}, Figure \ref{fig-sharp} and Figure \ref{fig-C1} for illustration.
In the case of sharp type, we can always shift $t_*$ to $0$ for convenience.

Our aim is to present a classification of those various types of wave solutions
with the admissible wave speeds depending on the time delay.
Throughout the paper we assume that
the death rate function $d(\cdot)$ satisfies
\begin{equation} \label{eq-d}
d\in C^2([0,+\infty)), \quad d(0)=0,~ d'(s)>0, ~d''(s)\ge0 \text{~for~} s>0,
\end{equation}
and the birth function $b$ satisfies
the following unimodality condition:
\begin{align} \nonumber
&b\in C^1(\mathbb{R}_{+};\mathbb{R}_{+})
\text{~ has only one positive local extremum point ~}
s=s_M
\\ \nonumber
&\text{(global maximum point) and~}
b(0)=0, b(\kappa)=d(\kappa), b'(0)>d'(0),
\\ \label{eq-UM}
&b'(\kappa)<d'(\kappa),
d(s)<b(s)\le b'(0)s \text{~for~} s\in(0,\kappa).
\end{align}

If $s_M\ge\kappa$, then $b$ is monotonically increasing on $[0,\kappa]$
and it is well known that the non-degenerate diffusion equation ($m=1$) admits
monotonically increasing wavefronts if and only if $c\ge c_*$ with
$c_*>0$ being the minimal wave speed determined by the characteristic equation
near the equilibrium $0$.
It is also shown in \cite{JDE18XU} that
the similar result holds for the degenerate diffusion equation ($m>1$)
except that the minimal wave speed is not determined by
the corresponding characteristic equation,
which indicates an essential difference between those two types of diffusion.
Henceforth, we may restrain ourselves to the case $s_M<\kappa$
in which $b$ is non-monotone in $[0,\kappa]$
and $b(s_M)>b(\kappa)=d(\kappa)$.
We set $M:=b(s_M)=\max b$, $\theta:=b(d^{-1}(M))$
and according to the monotone increasing of the death function $d$,
it holds $s_M<\kappa<d^{-1}(M)$.

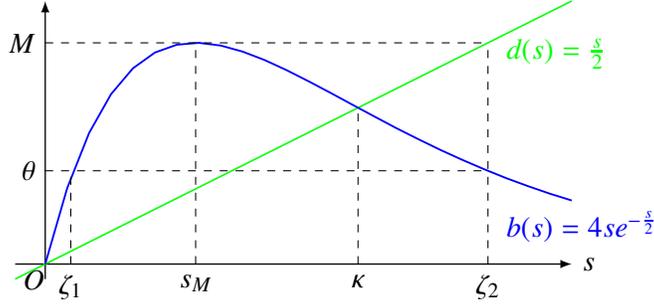
\begin{figure}[htb]
\begin{center}
\begin{tikzpicture}[scale=1.0,domain=0:7]
\def\x0{7} \def\y{3.5} \def\z{3} \def\p{4} \def\a{0.5}
\def\k{4.15888} \def\M{2.943} \def\z2{\M*2} \def\t{1.241} \def\w{0.34}
\draw[->,>=latex,line width=.5pt] (-0.4,0)--(\x0,0) node[right] {$s$};
\draw[->,>=latex,line width=.5pt] (0,-0.2)--(0,\y);
\node[below] at (-.15,.05) {$O$};
\draw[color=green,line width=.6pt] (-0.4,-0.2)--(\x0,\y);
\node[color=green,right] at (6,2.8) {$d(s)=\frac{s}{2}$};
\draw[color=blue,line width=.7pt] plot (\x,{\p*\x*exp(-\a*\x)});
\node[color=blue,right] at (6,0.5) {$b(s)=4se^{-\frac{s}{2}}$};
\draw[dashed,line width=.4pt] (\k,0)--(\k,{\k/2});
\node[below] at (\k,0) {$\kappa$};
\draw[dashed,line width=.4pt] (0,\M)--(\z2,\M);
\node[left] at (0,\M) {$M$};
\draw[dashed,line width=.4pt] (\z2,0)--(\z2,\M);
\node[below] at (\z2,0) {$\zeta_2$};
\draw[dashed,line width=.4pt] (0,\t)--(\z2,\t);
\node[left] at (0,\t) {$\theta$};
\draw[dashed,line width=.4pt] (\w,0)--(\w,\t);
\node[below] at (\w,0) {$\zeta_1$};
\draw[dashed,line width=.4pt] (2,0)--(2,\M);
\node[below] at (2,0) {$s_M$};
\end{tikzpicture}
\end{center}
\caption{The structure on functions $b(u)$ and $d(u)$.}
\label{fig-bd}
\end{figure}

The above unimodality condition \eqref{eq-UM} is stronger than the
following condition:
\begin{align} \nonumber
&{b : \mathbb{R}_{+} \rightarrow \mathbb{R}_{+}
\text { is continuous and such that, for some } 0<\zeta_{1}<\zeta_{2},}
\\ \nonumber
&{\text { 1. } b\left(\left[\zeta_{1}, \zeta_{2}\right]\right)
\subseteq\left[d(\zeta_{1}), d(\zeta_{2})\right]
\text { and } b\left(\left[0, \zeta_{1}\right]\right)
\subseteq\left[0, d(\zeta_{2})\right] ;}
\\ \nonumber
&{\text { 2. } \min _{s \in\left[\zeta_{1}, \zeta_{2}\right]} b(s)
=b\left(\zeta_{1}\right) ;}
\\ \nonumber
&{\text { 3. } b(s)>d(s) \text { for } s \in\left(0, \zeta_{1}\right]
\text { and } b \text { is differentiable at } 0, \text { with } b'(0)>d'(0);}
\\ \label{eq-zeta}
&{\text { 4. } \text{in}\left[0, \zeta_{2}\right],
\text { the equation } b(s)=d(s) \text { has exactly two solutions, } 0
\text { and } \kappa.}
\end{align}
Here we can take $\zeta_2=d^{-1}(M)=d^{-1}(\max b)$,
and $\zeta_1\in(0,s_M)$ such that $b(\zeta_1)=\theta$,
whose existence and uniqueness are ensured by the unimodality condition \eqref{eq-UM}
as shown in the illustrative Figure \ref{fig-bd}.

It is adapted for the case when the birth function $b(\cdot)$ satisfies
the following feedback condition:
\begin{equation}\label{eq-feedback}
(b(s)-\kappa)(s-\kappa)<0, \quad
s \in[d^{-1}(\theta), d^{-1}(M)] \backslash\{\kappa\}.
\end{equation}

Since the diffusion in \eqref{eq-tw} is degenerate for $\phi=0$,
and nonlinear with respect to $\phi$ near $\kappa$,
we define the following characteristic functions for $c>0$
near the two equilibria $0$ and $\kappa$ separately
\begin{equation} \label{eq-character}
\chi_{0}(\lambda):=b'(0)e^{-\lambda cr}-c\lambda-d'(0),
\qquad \lambda>0.
\end{equation}
and
\begin{equation}\label{eq-character-k}
\chi_{\kappa}(\lambda):= D m\kappa^{m-1}\lambda^{2}+b'(\kappa) e^{-\lambda cr}
-c \lambda-d'(\kappa), \qquad c>0.
\end{equation}

We see that $\chi_{0}(\lambda)=0$ has a unique positive real root
$\lambda_0$ for all $c>0$.
In fact, $\lambda_0=\frac{\omega_r}{c}$ such that $\omega_r\in(0,b'(0)-d'(0))$
is the unique solution of $b'(0)e^{-r\omega_r}=\omega_r+d'(0)$
since $b'(0)>d'(0)$.
However, the distribution of the roots of $\chi_{\kappa}(\lambda)$ is
much more complicated
and plays an essential role in determining the oscillatory property of the
semi-wavefronts.

For any given $m>1$, $D>0$ and $r\ge0$, we define the critical wave speed
$c_*(m,r,b,d)$ for the degenerate diffusion equation \eqref{eq-tw} as follows
\begin{align*}
c_*(m,r,b,d)
:=\inf\{c>0; \eqref{eq-tw} &\text{~admits semi-wavefronts (including wavefronts)}\}.
\end{align*}
Here we omit the dependence of the wave speed $c_*(m,r,b,d)$ on the parameter $D>0$
for simplicity since the dependence is trivial via a re-scaling method
such that the speed with $D>0$ is the speed with $D=1$ multiplied by $\sqrt{D}$.
This note is applicable for all the wave speeds in the rest of the paper.

Our main results are as follows.
First we state the existence and non-existence results of
wave solutions, including semi-wavefronts and wavefronts, sharp and smooth type.

\begin{theorem}[Existence of smooth waves] \label{th-exist}
For any $m>1$, $D>0$ and $r\ge0$, there exists a constant
$\hat c(m,r,b,d)>0$ depending on $m,r$ and the structure
of $b(\cdot),d(\cdot)$, such that for any $c>\hat c(m,r,b,d)$,
\eqref{eq-tw} admits smooth wave solutions $\phi(t)$
(semi-wavefronts or wavefronts, see Figure \ref{fig-smooth}) with
$$0<\zeta_1\le\liminf_{t\to+\infty}\phi(t)
\le\limsup_{t\to+\infty}\phi(t)\le\zeta_2,$$
and
$$
|\phi(t)-C_1e^{\lambda t}|\le C_2e^{\Lambda t}, \quad \text{~for any~} t<0,
$$
where $\lambda>0$ is the unique root of $\chi_0(\lambda)=0$ and
$\Lambda>\lambda$, $C_1,C_2>0$ are constants.
\end{theorem}

\begin{theorem}[Non-existence of waves] \label{th-non-exist}
For any $m>1$, $D>0$ and $r\ge0$, there exists a constant
$\dot c(m,r,b,d)>0$ depending on $m,r$ and the structure
of $b(\cdot),d(\cdot)$, such that, \eqref{eq-tw} admits no wave solution $\phi(t)$
(semi-wavefronts or wavefronts, sharp or smooth)
for any $c<\dot c(m,r,b,d)$.
Moreover,
$$\dot c(m,r,b,d)=\frac{\mu_0(m,b(\cdot),d(\cdot))+o(1)}{r},
\quad r\to+\infty,$$
where $\mu_0(m,b(\cdot),d(\cdot))>0$.
\end{theorem}

\begin{theorem}[Existence of sharp waves] \label{th-semifinite}
For any $m>1$, $D>0$ and $r\ge0$, there exists a constant
$c_0(m,r,b,d)>0$ depending on $m,r$ and the structure
of $b(\cdot),d(\cdot)$, such that for $c=c_0(m,r,b,d)$,
\eqref{eq-tw} admits sharp wave solutions $\phi(t)$
(semi-wavefronts or wavefronts, non-$C^1$ type (see Figure \ref{fig-sharp})
or $C^1$ type (see Figure \ref{fig-C1}))
with $\phi(t)\equiv0$ for $t\le0$,
$$0<\zeta_1\le\liminf_{t\to+\infty}\phi(t)
\le\limsup_{t\to+\infty}\phi(t)\le\zeta_2,$$
and
$$
|\phi(t)-C_1t_+^\lambda|\le C_2t_+^\Lambda, \quad \text{~for any~} t\in(0,1),
$$
where $t_+=\max\{t,0\}$, $\lambda=1/(m-1)$
and $\Lambda>\lambda$, $C_1,C_2>0$ are constants.
\end{theorem}

The sharp waves are classified into
$C^1$ type and non-$C^1$ type according to the degeneracy index $m$.

\begin{theorem}[Regularity of sharp waves] \label{th-sharp}
Assume that the conditions in Theorem \ref{th-semifinite} hold.
If $m\ge2$, then the sharp waves are of non-$C^1$ type
(as illustrated in Figure \ref{fig-sharp});
while if $1<m<2$, then the sharp traveling waves are of $C^1$ type
(as shown in Figure \ref{fig-C1}).
\end{theorem}

\begin{remark}
Roughly speaking, the degeneracy strengthens as $m>1$ increases
and the regularity of the case $m\ge2$ is weaker than that of $1<m<2$.
For the case $1<m<2$, the sharp traveling wave remains $C^1$ regularity
but not analytic.
\end{remark}

\begin{remark}
In the above theorems, we have introduced constants $\hat c(m,r,b,d)$,
$\dot c(m,r,b,d)$ and $c_0(m,r,b,d)$, with obviously,
$$
\dot c(m,r,b,d)\le c_*(m,r,b,d)\le \min\{\hat c(m,r,b,d),c_0(m,r,b,d)\},
$$
where $c_*(m,r,b,d)$ is the minimal wave speed,  or say critical wave speed.
We conjecture that the sharp type traveling wave is unique,  and the corresponding
wave speed
$$
c_*(m,r,b,d)=c_0(m,r,b,d).
$$
This is, the critical wave of the degenerate model is the unique sharp type
traveling wave,
and the speeds of smooth type wave solutions are greater than
the speed of sharp type wave solution.
Those conjectures are true for the case without time delay,
see for example \cite{Huang-Jin-Mei-Yin},
and for the case with time delay and quasi-monotonicity,
see our paper \cite{Our-Forthcoming}.
\end{remark}

\begin{remark}
The constants of speeds
$\hat c(m,r,b,d)$, $\dot c(m,r,b,d)$, $c_0(m,r,b,d)$ and $c_*(m,r,b,d)$
all are assumed to be dependent on the structure of functions $b(\cdot)$ and $d(\cdot)$.
It is well known that for the linear diffusion equation without time delay,
i.e., $m=1$, $r=0$, $c_*(1,0,b,d)=2\sqrt{D(b'(0)-d'(0))}$
provides that $b,d$ satisfy some concave structure.
Obviously, the critical wave speed of the linear diffusion equation is
totally determined by the linearization near zero.
However, the critical wave speed of the degenerate diffusion equation is
nonlinearly determined.
The wave front behaviors are controlled by the whole structure.
\end{remark}

Next we turn to the oscillating properties of wave solutions.

\begin{theorem}[Oscillating waves] \label{th-osc}
Assume that $m>1$, $r>0$, $b'(\kappa)<0$,
then there exists a number $c_\kappa=c_\kappa(m,r,b'(\kappa),d'(\kappa))\in(0,+\infty]$
such that the semi-wavefronts with speed $c>c_\kappa$
cannot be eventual monotone (i.e., they must be oscillating around $\kappa$,
convergent or divergent).
Moreover,
$$c_\kappa(m,r,b'(\kappa),d'(\kappa))=\frac{\mu_\kappa(m,b'(\kappa),d'(\kappa))+o(1)}{r},
\quad r\to+\infty,$$
where $\mu_\kappa(m,b'(\kappa),d'(\kappa)):=
\sqrt{\frac{2Dm\kappa^{m-1}\omega_\kappa}{b'(\kappa)}}e^{\frac{\omega_\kappa}{2}}$,
and $\omega_\kappa<-2$ is the unique negative root of
$2d'(\kappa)=b'(\kappa)e^{-\omega_\kappa}(2+\omega_\kappa)$.
\end{theorem}

\begin{theorem}[Non-decaying oscillating waves] \label{th-nondecay}
Assume that the function $b(\cdot)$ satisfies
the feedback condition \eqref{eq-feedback} and $b'(\kappa)<0$ and the time delay $r>0$,
then there exists a number $c^*=c^*(m,r,b'(\kappa),d'(\kappa))\in(0,+\infty]$ such that
the semi-wavefronts with speed $c>c^*$
have to develop non-decaying slow oscillations around $\kappa$.
Moreover, if $b'(\kappa)\ge -d'(\kappa)$,
then $c^*(m,r,b'(\kappa),d'(\kappa))=+\infty$ for large time delay $r$;
while if $b'(\kappa)<-d'(\kappa)$, then
$$
c^*(m,r,b'(\kappa),d'(\kappa))=
\frac{\mu^*(m,b'(\kappa),d'(\kappa))+o(1)}{r},
\quad r\to+\infty,
$$
where $\mu^*(m,b'(\kappa),d'(\kappa)):=
\pi\sqrt{\frac{Dm\kappa^{m-1}}{-b'(\kappa)-d'(\kappa)}}$.
\end{theorem}

\begin{remark}
For $b'(\kappa)\in[-d'(\kappa),0)$, we have $c^*(m,r,b'(\kappa),d'(\kappa))=+\infty$
and then $c^*(m,r,b'(\kappa),d'(\kappa))>c_\kappa(m,r,b'(\kappa),d'(\kappa))$
for large time delay.
For $b'(\kappa)\in(-\infty,-d'(\kappa))$, we also have
$c^*(m,r,b'(\kappa),d'(\kappa))>c_\kappa(m,r,b'(\kappa),d'(\kappa))$
for large time delay since
$$
\mu^*(m,b'(\kappa),d'(\kappa))=
\pi\sqrt{\frac{Dm\kappa^{m-1}}{-b'(\kappa)-d'(\kappa)}}
>
\sqrt{\frac{2Dm\kappa^{m-1}\omega_\kappa}{b'(\kappa)}}e^{\frac{\omega_\kappa}{2}}
=\mu_\kappa(m,b'(\kappa),d'(\kappa)),
$$
according the fact that $2|\omega_\kappa|e^{-|\omega_\kappa|}\le 2/e<\pi^2$
for all $\omega_\kappa$.
In fact, we show that
$c^*(m,r,b'(\kappa),d'(\kappa))\ge c_\kappa(m,r,b'(\kappa),d'(\kappa))$
for all cases in Lemma \ref{le-cstar}.
\end{remark}

\begin{remark}
In the above theorems, we investigate propagation dynamics of system \eqref{eq-main} without the
monotonicity assumption on the birth function $b(\cdot)$ for any large time delay $r$.
In the previous work \cite{Huang-Jin-Mei-Yin}, the authors proved the existence of traveling waves solutions for small time delay due to the limitation of perturbation method.
\end{remark}

To conclude, the time delay $r$ and  the degenerate diffusion have a strong influence on the geometry of wave profiles.
Here, we depict the shape of the traveling waves characterized by the wave speed $c$
and time delay $r$.
From Theorem \ref{th-exist}, we know that when the wave speed $c>\hat c$,
there exist smooth traveling wave solutions.
Theorem \ref{th-semifinite} implies there exists a sharp traveling wave with the wave speed $c_0$.
After investigate the geometry of leading edge,
it is naturally to consider the convergence of the semi-wavefronts.
So we have Theorem \ref{th-osc} and Theorem \ref{th-nondecay} indicating the oscillating
properties for both the sharp type and smooth type traveling waves.

\begin{figure}[htb]
\begin{center}
\begin{tikzpicture}[scale=1.4,domain=0:8]
\def\axl{0} \def\axr{8} \def\ay{4}
\draw[->,>=latex,line width=.5pt] (\axl,0)--(\axr,0);
\node at (7,-.4) {time-delay $r$};
\draw[->,>=latex,line width=.5pt] (0,0)--(0,\ay);
\node at (1,3.8) {wave speed $c$};
\draw[loosely dotted,color=blue,line width=.9pt]
    plot[domain=1:8, samples=144, smooth] (\x,{3.5/(\x)^0.3});
\node at (3,2.8) {$c^*$};
\draw[densely dotted,color=blue,line width=.8pt]
    plot[domain=0.5:8, samples=144, smooth] (\x,{3/(\x/0.5)^0.5});
\node at (2,1.68) {$c_\kappa$};
\draw[dashed,color=blue,line width=.7pt]
    plot[domain=0:8, samples=144, smooth] (\x,{4/(2+\x)});
\node at (1,1.54) {$\hat c$};
\draw[color=blue,line width=.7pt]
    plot[domain=0:8, samples=144, smooth] (\x,{2/(2+\x)});
\node at (0.5,0.6) {$c_0$};
\end{tikzpicture}
\end{center}
\caption{Different types of traveling waves for the degenerate diffusion equation
with time delay \eqref{eq-main}
correspond to time delay $r$ and wave speed $c$:
the case that the curve $c_0(m,r,b,d)$ never intersects with
the curves $c_\kappa(m,r,b'(\kappa),d'(\kappa))$ and $c^*(m,r,b'(\kappa),d'(\kappa))$.
}
\label{fig-cr}
\end{figure}
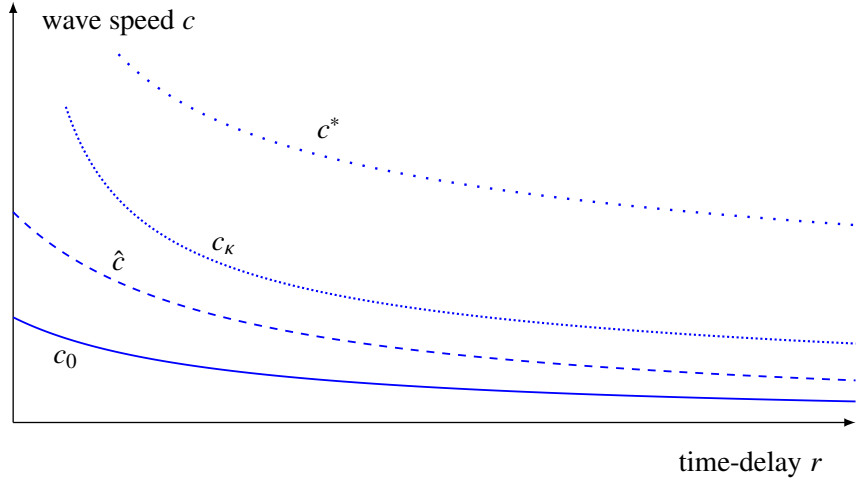

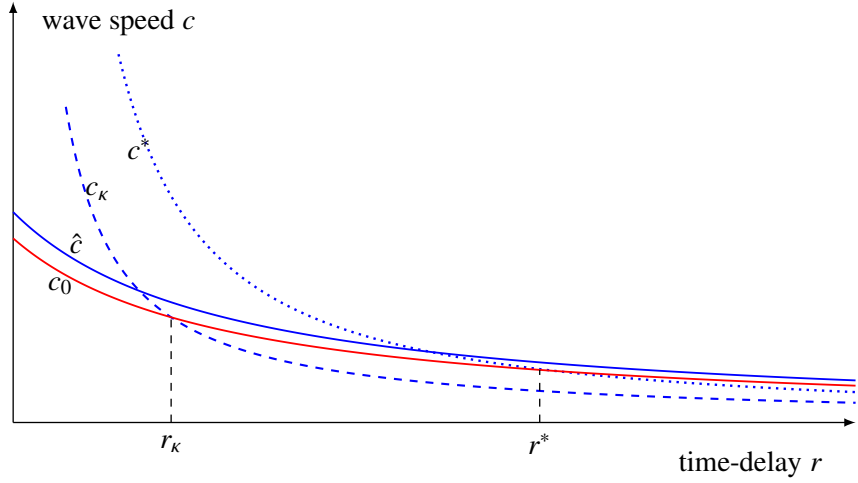
\begin{figure}[htb]
\begin{center}
\begin{tikzpicture}[scale=1.4,domain=0:8]
\def\axl{0} \def\axr{8} \def\ay{4}
\draw[->,>=latex,line width=.5pt] (\axl,0)--(\axr,0);
\node at (7,-.4) {time-delay $r$};
\draw[->,>=latex,line width=.5pt] (0,0)--(0,\ay);
\node at (1,3.8) {wave speed $c$};
\draw[dotted,color=blue,line width=.9pt]
    plot[domain=1:8, samples=144, smooth] (\x,{3.5/(\x)^1.2});
\node at (1.2,2.6) {$c^*$};
\draw[dashed,color=blue,line width=.8pt]
    plot[domain=0.5:8, samples=144, smooth] (\x,{3/(\x/0.5)^1});
\node at (0.8,2.2) {$c_\kappa$};
\draw[color=blue,line width=.7pt]
    plot[domain=0:8, samples=144, smooth] (\x,{4/(2+\x)});
\node at (0.6,1.7) {$\hat c$};
\draw[color=red,line width=.7pt]
    plot[domain=0:8, samples=144, smooth] (\x,{3.5/(2+\x)});
\node at (0.45,1.3) {$c_0$};
\draw[dashed,line width=.5] (1.5,0)--(1.5,1) node at (1.5,-.2) {$r_\kappa$};
\draw[dashed,line width=.5] (5,0)--(5,0.5) node at (5,-0.2) {$r^*$};
\end{tikzpicture}
\end{center}
\caption{Different types of traveling waves for the degenerate diffusion
equation with time delay \eqref{eq-main}
correspond to time delay $r$ and wave speed $c$:
the case that the curve $c_0(m,r,b,d)$ intersects with
the curves $c_\kappa(m,r,b'(\kappa),d'(\kappa))$ and $c^*(m,r,b'(\kappa),d'(\kappa))$
at $r_\kappa$ and $r^*$ respectively.}
\label{fig-cr-2}
\end{figure}

According to the above theorems, Figure \ref{fig-cr}
and Figure \ref{fig-cr-2} illustrate two main possible sketches of the
corresponding wave behaviors varying with the traveling wave speed $c$ and time delay $r$.
The critical lines of the wave speeds depend on the time delay and divide the $(r, c)$
plane into several parts relating to different wave behaviors.
The slopes and structures of these curves depend on the functions $b(\cdot)$ and $d(\cdot)$.
It is worth to mention that there exist sharp-oscillating waves
for some proper parameters, which is different from the former literatures
(see Figure \ref{fig-sharp} and Figure \ref{fig-C1}).

If the curve $c_0(m,r,b,d)$ never intersects with
the curves $c_\kappa(m,r,b'(\kappa),d'(\kappa))$ and $c^*(m,r,b'(\kappa),d'(\kappa))$
as illustrated in Figure \ref{fig-cr},
we have the following different types of waves:
the curve $c_0$ is the wave speed of sharp type traveling waves;
the waves with the parameters $(r,c)$ above the curve $\hat c$ are positive and smooth
and the types (A1), (A2) and (A3) in Figure \ref{fig-smooth} are possible;
the waves with $(r,c)$ above the curve $c_\kappa$ are oscillatory;
the waves with $(r,c)$ above the curve $c^*$ are non-decaying oscillatory.
If $b'(\kappa)\ge0$, then $c_\kappa(m,r,b'(\kappa),d'(\kappa))=+\infty$
and $c^*(m,r,b'(\kappa),d'(\kappa))=+\infty$, and
the curve $c_0(m,r,b,d)$ never intersects with the curves
$c_\kappa(m,r,b'(\kappa),d'(\kappa))$ and $c^*(m,r,b'(\kappa),d'(\kappa))$.
Actually, for the monotonically increasing function $b(\cdot)$,
the traveling waves are monotone.

Wave dynamics are rather complicated
when the curve $c_0(m,r,b,d)$ or $\hat c(m,r,b,d)$ intersects with the curve
$c_\kappa(m,r,b'(\kappa),d'(\kappa))$ or $c^*(m,r,b'(\kappa),d'(\kappa))$.
It raises the possibility of nine types of traveling waves as shown
in Figure \ref{fig-smooth}-\ref{fig-C1}.
Figure \ref{fig-cr-2} shows the
case that the curve $c_0(m,r,b,d)$ intersects with
the curves $c_\kappa(m,r,b'(\kappa),d'(\kappa))$ and $c^*(m,r,b'(\kappa),d'(\kappa))$
at $r_\kappa$ and $r^*$ respectively.
In this situation, many types of waves occur depending on the wave speed $c$,
the time delay $r$ and the degeneracy $m$ as follows:
\\ \indent
(i) if the degeneracy is strong with $m\ge2$, then along the curve
$c_0(m,r,b,d)$, the non-$C^1$ sharp type wave is monotone (B1) for small time delay
or non-monotone (B2) if $r>r_\kappa$ or non-decaying oscillatory (B3)
if $r>r^*$;
\\ \indent
(ii) if the degeneracy is weak with $1<m<2$, then along the curve
$c_0(m,r,b,d)$, the sharp waves are $C^1$ type, that is,
(B1), (B2), (B3) are replaced by (C1), (C2) and (C3);
\\ \indent
(iii) the waves with the parameters $(r,c)$ above the curve $\hat c(m,r,b,d)$
are positive and smooth, that is, (A1), (A2), (A3) are possible if the time delay is small;
\\ \indent
(iv) after the curve $\hat c(m,r,b,d)$ intersects
with $c_\kappa(m,r,b'(\kappa),d'(\kappa))$ (it happens if $\dot c(m,r,b,d)$ intersects
with $c_\kappa(m,r,b'(\kappa),d'(\kappa))$ since $\hat c\ge \dot c$),
the monotone waves are impossible, that is, only (A2) and (A3) of smooth type exist;
\\ \indent
(v) after the curve $\hat c(m,r,b,d)$ intersects
with $c^*(m,r,b'(\kappa),d'(\kappa))$ (it happens if $\dot c(m,r,b,d)$ intersects
with $c^*(m,r,b'(\kappa),d'(\kappa))$),
the smooth wave has to develop non-decaying oscillations, that is,
only (A3) of the smooth type exists.

\begin{remark}
We note that for functions $b(\cdot)$ and $d(\cdot)$ with some structure condition,
the curves $c_0(m,r,b,d)$ and $\hat c(m,r,b,d)$ intersect with
$c_\kappa(m,r,b'(\kappa),d'(\kappa))$ and $c^*(m,r,b'(\kappa),d'(\kappa))$,
and then the various types of waves can happen.
In fact, $\mu_\kappa(m,b'(\kappa),d'(\kappa))$ and $\mu^*(m,b'(\kappa),d'(\kappa))$
are constants only depending on the local property of $b'(\kappa)$
and converge to zero as $b'(\kappa)\to-\infty$.
From the proof of Theorem \ref{th-non-exist}, we see that
$\mu_0(m,b,d)$ depends on the structure of $b(\cdot)$ and $d(\cdot)$
within $(0,\zeta_1)$, where $\zeta_1$ is determined
by the whole structure of $b(\cdot)$ and $d(\cdot)$
as shown in Figure \ref{fig-bd}.
The local variation of $b'(\kappa)$ has minor effect on $\mu_0(m,b,d)$
(if the change of $b'(\kappa)$ has no effect on $\zeta_1$, then $\mu_0(m,b,d)$ is fixed).
Hence, for functions $b(\cdot)$ with appropriate structure near $0$ and
suitable large $-b'(\kappa)$, there holds
$$
0<\mu_\kappa(m,b'(\kappa),d'(\kappa))<\mu^*(m,b'(\kappa),d'(\kappa))
<\mu_0(m,b(\cdot),d(\cdot)),
$$
and further
$$
0<c_\kappa(m,r,b'(\kappa),d'(\kappa))
<c^*(m,r,b'(\kappa),d'(\kappa))
<\dot c(m,r,b,d)\le c_0(m,r,b,d),
$$
for large time delay according to the asymptotic behavior in
Theorem \ref{th-non-exist}, Theorem \ref{th-osc}
and Theorem \ref{th-nondecay}.
\end{remark}

\section{Existence of traveling wave solutions}
In this section,
we employ the Schauder's Fixed Points Theorem to show the existence of
monotone and non-monotone traveling wave solutions.
Compared with the linear diffusion case ($m=1$),
both the comparison principle and the solvability of degenerate
elliptic problem ($m>1$) are not obvious.
The solvability of linear diffusion case can be shown by
writing the explicit expression by applying the variation of constants formula.
We can not expect such kind of expressions due to the degenerate diffusion.

Here we recall the comparison principle of degenerate diffusion
equation on unbounded domain proved in \cite{JDE18XU}.

\begin{lemma}[Comparison Principle, \cite{JDE18XU}] \label{le-comparison}
Let $\phi_1,\phi_2\in C(\mathbb R;\mathbb R)$
such that for $i=1,2$,
$0\le\phi_i\in L^\infty(\mathbb R)$, $\phi_i^m\in W_{\text{loc}}^{1,2}$,
$\phi_1(t)>0$ for all $t\in\mathbb R$,
$\phi_i(t)$ is increasing for $t\le t_0$ with some fixed $t_0\in\mathbb R$,
$\liminf_{t\to\pm\infty}(\phi_1(t)-\phi_2(t))\ge0$,
$\liminf_{t\to+\infty}\phi_1(t)>0$
and $\phi_i$ satisfies the following inequality
$$
c\phi_1'(t)-D(\phi_1^m(t))''+d(\phi_1(t))
\ge
c\phi_2'(t)-D(\phi_2^m(t))''+d(\phi_2(t))
$$
in the sense of distributions.
Then $\phi_1(t)\ge\phi_2(t)$ for all $t\in\mathbb R$.
\end{lemma}

We also need the following solvability and monotonicity of degenerate equations
on unbounded domain.

\begin{lemma} \label{le-solv}
Assume that $0\le\psi(t)\in L^\infty(\mathbb R)\cap C(\mathbb R)$,
$\psi$ is monotonically increasing on $(-\infty,t_0]$ for some $t_0\in\mathbb R$,
and $\psi(t)\ge\psi(t_0)>0$ for all $t>t_0$,
then the following degenerate elliptic equation
\begin{equation} \label{eq-solv}
\begin{cases}
c\phi'(t)-D(\phi^m(t))''+d(\phi(t))=\psi(t), \quad t\in\mathbb R, \\
\displaystyle
\lim_{t\to-\infty}\phi(t)=0, \\
\displaystyle
0<d^{-1}(\liminf_{t\to+\infty}\psi(t))\le
\liminf_{t\to+\infty}\phi(t)\le\limsup_{t\to+\infty}\phi(t)
\le d^{-1}(\limsup_{t\to+\infty}\psi(t))<+\infty,
\end{cases}
\end{equation}
admits at least one solution $\phi(t)$ such that
$0\le\phi(t)\in L^\infty(\mathbb R)$, $\phi$ is monotonically increasing
on $(-\infty,t_0]$,
and $\phi(t)\ge\phi(t_0)>0$ for all $t>t_0$.
\end{lemma}
{\it \bfseries Proof.}
This proof is similar to that of Lemma 3.5 in \cite{JDE18XU}.
Consider the following regularized problem for any $A>\max\{1,t_0\}$ with $-A<t_0$
\begin{equation} \label{eq-zregularized}
\begin{cases}
\displaystyle
c\phi'(t)=D\big(m(|\phi(t)|^2+1/A)^{(m-1)/2}\phi'(t)\big)'-d(\phi(t))
+\psi(t), \quad t\in(-A,A), \\
\displaystyle
\phi(-A)=d^{-1}(\psi(-A)), \quad \phi(A)=d^{-1}(\psi(A)).
\end{cases}
\end{equation}
The unique existence of solution to \eqref{eq-zregularized} is trivial.
The solution is denoted by $\phi_A$.
We note that $d(s)$ is monotonically increasing and
$\psi(t)\ge\psi(-A)$ for all $t\ge -A$ since $-A<t_0$ and
$\psi(t)$ is increasing on $(-\infty,t_0)$.
Comparison principle of elliptic equation shows that
$$0<d^{-1}(\psi(-A))\le\phi_A(t)\le d^{-1}(\sup \psi), \quad t\in(-A,A).$$
In fact, if this is not true, we argue by contradiction.
If there exists $t_0\in(-A,A)$ such that $\phi_A(t_0)<d^{-1}(\psi(-A))$,
then the minimum of $\phi_A(t)$ on $[-A,A]$ is less than $d^{-1}(\psi(-A))$
and is attained at some inner point $t^*\in(-A,A)$
since at the endpoints $\phi_A(\pm A)\ge d^{-1}(\psi(-A))$.
At this point $t^*$, $\phi_A'(t^*)=0$, $\phi_A''(t^*)\ge0$,
and by \eqref{eq-zregularized}
$$\psi(t^*)=c\phi_A'(t^*)-D\big(m(|\phi_A(t^*)|^2+1/A)^{(m-1)/2}\phi_A'(t^*)\big)'
+d(\phi_A(t^*))<\psi(-A),$$
which contradicts to the fact $\psi(t)\ge \psi(-A)$ for all $t\in[-A,A]$.
The proof of $\phi_A(t)\le d^{-1}(\sup\psi)$ is similar.

We assert that $\phi_A'(t)\ge0$ for $t\in[-A,t_0]$.
Otherwise, there exists a $t_*\in(-A,t_0)$ such that $\phi_A'(t_*)<0$.
Let $(t_1,t_2)$ be the maximal interval such that $t_*\in(t_1,t_2)$
and $\phi_A'(t)<0$ for $t\in(t_1,t_2)$.
We note that $\phi_A(t)$ attains its minimum at $-A$,
which implies $\phi_A'(-A)\ge0$.
Thus, $t_1\in[-A,t_*)$, $\phi_A'(t_1)=0$,
$$(m(|\phi_A(t)|^2+1/A)^{(m-1)/2}\phi_A'(t))'\big|_{t=t_1}\le0, $$
and
\begin{align*}
\psi(t_1)=c\phi_A'(t_1)
-D\big(m(|\phi_A(t_1)|^2+1/A)^{(m-1)/2}\phi_A'(t_1)\big)'+d(\phi_A(t_1))
\ge d(\phi_A(t_1)),
\end{align*}
which shows
$$
\phi_A(t_1)\le d^{-1}(\psi(t_1))\le d^{-1}(\psi(t_*))\le d^{-1}(\psi(t_0))
\le d^{-1}(\psi(A))=\phi_A(A)
$$
as $t_1\le t_*<t_0<A$, $\psi(t)$ is
increasing on $(-\infty,t_0]$ and $\psi(t)>\psi(t_0)>0$ for all $t>t_0$.
Therefore, $\phi_A(A)=d^{-1}(\psi(A))\ge \phi_A(t_1)$ and
$\phi_A(t)$ cannot always decreasing on the whole $(t_1,A)$.
Then $t_2<A$ and $\phi_A'(t_2)=0$, $\phi_A(t_1)>\phi_A(t_2)$,
$$(m(|\phi_A(t)|^2+1/A)^{(m-1)/2}\phi_A'(t))'\big|_{t=t_2}\ge0,$$
and
\begin{align*}
\psi(t_1)&=c\phi_A'(t_1)
-D\big(m(|\phi_A(t_1)|^2+1/A)^{(m-1)/2}\phi_A'(t_1)\big)'+d(\phi_A(t_1))\\
&>c\phi_A'(t_2)
-D\big(m(|\phi_A(t_2)|^2+1/A)^{(m-1)/2}\phi_A'(t_2)\big)'+d(\phi_A(t_2))\\
&=\psi(t_2), \quad t_1<t_2,
\end{align*}
which contradicts to the monotonically increasing of $\psi$ on $(-\infty,t_0)$
and $\psi(t)\ge\psi(t_0)>0$ for all $t>t_0$.

Next, we show that $\phi_A(t)\ge\phi_A(t_0)>0$ for all $t>t_0$.
Otherwise, there exists a number $t_1\in(t_0,A)$ such that $\phi_A(t_1)<\phi_A(t_0)$.
Noticing that $\phi_A(t)$ is increasing on $(-A,t_0)$,
we see that there exists a maximum point $t^*\in[t_0,t_1)$.
Similar to the above analysis at this point $t^*$, we find that
$$\phi_A(A)\ge\phi_A(t^*)\ge \phi_A(t_0)>\phi_A(t_1)$$
and $\phi_A(t)$ cannot decreasing on the whole $(t^*,A)$.
Then their exist $t_a\in(t^*,t_1)$ and $t_b\in(t_1,A)$ such that
$\phi_A(t_a)=\phi_A(t_b)=\phi_A(t_0)$
and $\phi_A(t)$ satisfies
\begin{equation*}
\begin{cases}
\displaystyle
c\phi'(t)=D\big(m(|\phi(t)|^2+1/A)^{(m-1)/2}\phi'(t)\big)'-d(\phi(t))
+\psi(t), \quad t\in(t_a,t_b), \\
\displaystyle
\phi(t_a)=\phi_A(t_0), \quad \phi(t_b)=\phi_A(t_0).
\end{cases}
\end{equation*}
Applying the maximum principle of elliptic equations with $\psi(t)\ge \psi(t_0)$
for all $t>t_0$, we find that
$\phi_A(t)\ge \phi_A(t_0)$ for $t\in(t_a,t_b)$,
which contradicts to $t_1\in(t_a,t_b)$ and $\phi_A(t_1)<\phi_A(t_0)$.

For any $1<B<A$, let $\eta(t)$ be the cut-off function such that
$0\le\eta(t)\le1$, $\eta\in C_0^2((-B,B))$, $|\eta'(t)|\le2$ for $t\in(-B,B)$,
$\eta(t)=1$ for $t\in(-B+1,B-1)$.
Multiply \eqref{eq-zregularized} by $\eta^2(t)\phi_A(t)$ and integrate over $(-A,A)$,
we have
\begin{eqnarray*}
&& \int_{-A}^A c\eta^2 \phi_A(t)\phi_A'(t)dt
+\int_{-A}^ADm\eta^2(|\phi_A(t)|^2+1/A)^{(m-1)/2}|\phi_A'(t)|^2dt\\
&& \ \ \ +\int_{-A}^A\eta^2 d(\phi_A(t))\phi_A(t)dt\\
&&\le \int_{-A}^A2Dm\eta(|\phi_A(t)|^2+1/A)^{(m-1)/2}\phi_A(t)\phi_A'(t)|\eta'(t)|dt
+\int_{-A}^A\eta^2\phi_A(t)\psi(t)dt\\
&&\le\frac{1}{2}\int_{-A}^ADm\eta^2(|\phi_A(t)|^2+1/A)^{(m-1)/2}|\phi_A'(t)|^2dt\\
&&\ \ \ +\int_{-A}^A2Dm(|\phi_A(t)|^2+1/A)^{(m-1)/2}|\phi_A(t)|^2|\eta'(t)|^2dt
+2d^{-1}(\sup\psi)\sup\psi B.
\end{eqnarray*}
Therefore,
\begin{align*}
&\frac{1}{2}\int_{-B+1}^{B-1}Dm(|\phi_A(t)|^2+1/A)^{(m-1)/2}|\phi_A'(t)|^2dt
+\int_{-B+1}^{B-1}d(\phi_A(t))\phi_A(t)dt\\
&\le\int_{-B}^{-B+1}+\int_{B-1}^{B}
2Dm(|\phi_A(t)|^2+1/A)^{(m-1)/2}|\phi_A(t)|^2|\eta'(t)|^2dt
+2d^{-1}(\sup\psi)\sup\psi B\\
&\le 16Dm((\sup\psi)^2+1)^{(m-1)/2}(\sup\psi)^2+2d^{-1}(\sup\psi)\sup\psi B.
\end{align*}
It follows that $\|\phi_A^m\|_{W^{1,2}(-B+1,B-1)}$ is uniformly bounded and independent of $A$.
We note that the embedding $W^{1,2}(-B+1,B-1)$ into $C^{\gamma}([-B+1,B-1])$
with $\gamma\in(0,\frac{1}{2})$ is compact,
and $\phi_A^m\in C^{\gamma}([-B+1,B-1])$ implies $\phi_A\in C^{\gamma/m}([-B+1,B-1])$.
There exist a subsequence of $\{\phi_A(t)\}_{A>1}$
denoted by $\{\phi_{A_n}(t)\}_{n\in\mathbb N}$ and
a function $\phi(t)\in C^{\gamma/m}(\mathbb R)$ such that
$\phi^m\in W_\mathrm{loc}^{1,2}(\mathbb R)$, $0\le\phi\le K$,
and $\phi_{A_n}(t)$ uniformly converges to $\phi(t)$ on any compact interval,
$\phi_{A_n}^m(t)$ weakly converges to $\phi^m(t)$ in $W_\mathrm{loc}^{1,2}(\mathbb R)$.
Since each $\phi_{A_n}(t)$ is monotonically increasing on $(-\infty,t_0)$, we see that
$\phi(t)$ is also increasing on $(-\infty,t_0)$.
We can verify that $\phi(t)$ is a solution of \eqref{eq-solv}.
$\hfill\Box$

The following uniform permanence property is
similar to the linear diffusion case ($m=1$) in
\cite{TrofimchukJDE, Trofimchuk-Admissible}.
Their proof was based on the variation of constants formula for
semilinear differential equations.
Here we adopt an alternative proof applicable for
quasi-linear diffusion equations ($m>1$).

\begin{lemma}[Uniform permanence property] \label{le-permanence}
Assume the unimodality condition \eqref{eq-UM}
or its consequence \eqref{eq-zeta} with additionally
supposing that $\sup_{s\ge0}b(s)\le M$.
Then any non-trivial wave solution $\phi(t)$
of \eqref{eq-main} satisfies
$$
0<\zeta_1\le\liminf_{t\to+\infty}\phi(t)
\le\limsup_{t\to+\infty}\phi(t)\le\zeta_2<+\infty.
$$
\end{lemma}
{\it\bfseries Proof.}
We first prove that $\limsup_{t\to+\infty}\phi(t)\le\zeta_2$.
We proceed by contrary, supposing that there exists a $t_*\in\mathbb R$ such that
$\phi(t_*)>\zeta_2$.
Let $(t_1,t_2)$ be the maximal interval such that $t_*\in(t_1,t_2)$ and
$\phi(t)>\zeta_2$ in $(t_1,t_2)$, i.e.,
$(t_1,t_2)$ is the connected component containing $t_*$
of the set that $\phi(t)>\zeta_2$.
Since $\lim_{t\to-\infty}\phi(t)=0<\zeta_2$, we see that $t_1\in(-\infty,t_*)$.
If there is no local maximum point of $\phi(t)$ in $(t_1,t_2)$, then $t_2=+\infty$ and
$\phi(t)$ is monotonically increasing and converges to some equilibrium
greater than $\zeta_2$ as $t$ tends to positive infinity,
which is impossible since the only positive equilibrium is $\kappa<\zeta_2$.
Now let $t_0\in (t_1,t_2)$ be one of the local maximum points.
We have $\phi(t_0)\ge\phi(t_*)>\zeta_2$, $\phi'(t_0)=0$,
$(\phi^m(t))''|_{t=t_0}\le0$ as $t_0$ is also a maximum point of $\phi^m(t)$,
and at this point $t_0$
\begin{align*}
c\phi'(t)-D(\phi^m(t))''&+d(\phi(t))-b(\phi(t-cr))\\
&\ge d(\phi(t))-b(\phi(t-cr))
> d(\zeta_2)-M=0,
\end{align*}
which contradicts to the equation \eqref{eq-tw}.
Therefore, we proved that $\phi(t)\le\zeta_2$ for all $t\in\mathbb R$.

We next prove that $\phi(t)>0$ for $t\in\mathbb R$ unless $\phi(s)\equiv0$ for
all $s\le t$, which is in fact possible for the sharp type wave solution.
Suppose that there exists a $t_*$ such that $\phi(t_*)=0$ and
$\phi(s)\not\equiv0$ for $s\le t_*$.
Here at $t_*$, we have $\phi'(t)=0$, $(\phi^m(t))''\ge0$ and
\begin{align*}
b(\phi(t-cr))=c\phi'(t)-D(\phi^m(t))''+d(\phi(t))\le0,
\end{align*}
which means $\phi(t_*-cr)=0$
and $\phi(t_*-jcr)=0$ for all $j\in \mathbb Z^+$ by induction.
Supposing that $t_0$ is the boundary of the leading edge of $\phi(t)$
(see Definition \ref{de-tw}) and in this case $t_0<t_*<+\infty$, $\phi(t_0)>\kappa$,
$\phi'(t_0)=0$, $(\phi^m(t))''|_{t=t_0}\le0$, then we have at $t_0$
\begin{align*}
b(\phi(t-cr))=c\phi'(t)-D(\phi^m(t))''+d(\phi(t))> d(\kappa)>0.
\end{align*}
It follows that $\phi(t_0-cr)>0$ and $\phi(t)>0$ for $t\in(t_0-cr,t_0)$ since
$\phi$ is monotonically increasing in $(-\infty,t_0)$.
That is, we find an interval in $(-\infty,t_*)$
longer than $cr$ such that $\phi(t)$ has no zero point,
which contradicts to $\phi(t_*-jcr)=0$ for all $j\in \mathbb Z^+$.
We conclude that $\phi(t)>0$ for all $t>t_0$.

We finally prove that $\zeta_1\le\liminf_{t\to+\infty}\phi(t)$.
Assuming that $\liminf_{t\to+\infty}\phi(t)<\zeta_1$, then
there exists a sequence $\{t_n\}_{n=1}^{\infty}$ such that
$t_n$ tends to positive infinity and $\phi(t_n)<\zeta_1$.
Let $A=\{t>t_1; \phi(t)<\zeta_1\}$.
We denote the set of all the local minimum points
of $\phi(t)$ in $A$ by $A_\text{min}$.
We divide the following proof into two parts.

(i) If $A_\text{min}$ is empty or bounded to the upwards,
then $\phi(t)$ is eventually monotone and
converges to some equilibrium in $[0,\zeta_1]$,
which can only be $0$.
Therefore, $\phi(t)$ is monotonically decreasing on $[t_*,+\infty)$
and $\phi(t)\in[0,\varepsilon]$ for some sufficiently large $t_*$,
where $\varepsilon\in(0,\zeta_1)$ such that
$$\frac{b(s)-d(s)}{s}\ge \frac{b'(0)-d'(0)}{2}, \quad \forall s\in(0,\varepsilon)$$
since the limit of the left hand side is $b'(0)-d'(0)>0$ as $s$ tends to zero.
Now we have for $t>t_*+cr$, $\phi'(t)\le0$, $\phi(t-cr)\ge\phi(t)$ and
$$b(\phi(t-cr))\ge b(\phi(t))\ge d(\phi(t))+\frac{b'(0)-d'(0)}{2}\phi(t).$$
Here we have used the monotonicity of $b(s)$ on $[0,\varepsilon]$
since $b'(0)>0$ and we may take $\varepsilon$ even smaller if necessary.
Then
\begin{equation} \label{eq-zsubode}
D(\phi^m(t))''=c\phi'(t)+d(\phi(t))-b(\phi(t-cr))
\le -\frac{b'(0)-d'(0)}{2}\phi(t),
\end{equation}
which must decay to zero at some finite point $\hat t\in (t_*,+\infty)$
and $\phi(t)>0$ for $t<\hat t$
according to the phase plane analysis for this sublinear ordinary differential
equation \eqref{eq-zsubode}.
At this point $\hat t$, we also have
\begin{align*}
b(\phi(t-cr))=c\phi'(t)-D(\phi^m(t))''+d(\phi(t))=0,
\end{align*}
which contradicts to $\phi(\hat t-cr)>0$.

(ii) If $A_\text{min}$ is unbounded to the upwards.
Suppose that $t_0$ is the boundary of the leading edge of $\phi(t)$
(see Definition \ref{de-tw}) and in this case
$t_0<t_1\le \inf A\le \inf A_\text{min}<+\infty$, $\phi(t_0)>\kappa$,
($t_1>t_0$ is trivial as we can modify the sequence $\{t_n\}_{n=1}^\infty$).
We have already proved that $\phi(t)>0$ for $t>t_0$, and the local regularity
of non-degenerate diffusion equation \eqref{eq-tw} shows that
there is no bounded accumulation point of $A_\text{min}$.
For any $t_*\in A_\text{min}$, we find that $\phi'(t_*)=0$,
$(\phi^m(t))''|_{t=t_*}\ge0$ and
\begin{equation} \label{eq-zperm}
b(\phi(t-cr))=c\phi'(t)-D(\phi^m(t))''+d(\phi(t))\le d(\phi(t)) \text{~~at~} t_*.
\end{equation}
According to the structure assumption on $b(s)$ and $d(s)$, we can choose
positive constants $k_1>0$ and $k_2\in(0,1)$ such that
$b(s)-d(s)\ge k_1 s$ and $b(\tilde s)\ge d(s)$ for all $s\in(0,\zeta_1)$
and $\tilde s\in[k_2s,s]$.
We deduce from \eqref{eq-zperm} that
$\phi(t_*-cr)<k_2 \phi(t_*)<\phi(t_*)$.
Noticing that $t_*$ is a local minimum point, we see that
$t_*-cr<t_1$ or
there exists another local minimum point $\tilde t_*<t_*$
and $\tilde t_*\in A_\text{min}$ such that
$\phi(\tilde t_*)\le \phi(t_*-cr)<k_2 \phi(t_*)$,
which is denoted by $\tilde t_*=F(t_*)$ for convenience.
Furthermore, if $F(t_*^1)=F(t_*^2)$ for two different minimum
points $t_*^1,t_*^2\in A_\text{min}$ with $t_*^1<t_*^2$,
then $0<t_*^2-t_*^1<cr$ as $t_*^2-cr<t_*^1$,
otherwise, $F(t_*^2)\ge t_*^1>F(t_*^1)$, a contradiction.
Therefore, we can choose a subsequence
$\{s_n\}_{s=1}^\infty\subset A_\text{min}$ such that
$\phi(s_{n+1})\ge \phi(s_n)/k_2$ for all $n\in\mathbb Z^+$,
which contradicts to $k_2\in(0,1)$ and $\phi(t)<\zeta_1$ for all $t\in A$.
The proof is completed.
$\hfill\Box$

The existence of traveling waves is deduced by Schauder fixed point theorem
on an appropriate profile set $\Gamma_\epsilon$ constructed with
upper and lower profiles $\phi^*$ and $\phi_\epsilon$
for two auxiliary problems,
where $\phi^*,\phi_\epsilon$ will be specified in the following.
We follow the main lines of Theorem 1.1 in \cite{MaJDE}
and construct two auxiliary reaction diffusion equations with quasi-monotonicity.
Since $b(\zeta_1)>0$, there is a small $\epsilon_0\in(0,\zeta_1)$
such that $b(\zeta_1-\epsilon)>0$ for every $\epsilon\in[0,\epsilon_0]$.
If $b(s)$ satisfies \eqref{eq-zeta}, for any $\epsilon\in(0,\epsilon_0)$,
define two continuous functions as follows
$$
b^{*}(u)=
\begin{cases}
\min \left\{b'(0) u,M\right\}, \quad u\in[0,\zeta_2],\\
\max \left\{M, b(u)\right\}, \quad u>\zeta_2,
\end{cases}
$$
and
$$
b_{\epsilon}(u)=
\begin{cases}
\displaystyle
\inf_{\eta \in\left[u, \zeta_2\right]}\left\{b(\eta), d\left(\zeta_1-\epsilon\right)\right\}, \quad u\in[0,\zeta_2],\\
\min \left\{b(u), d\left(\zeta_1-\epsilon\right)\right\}, \quad u>\zeta_2.
\end{cases}
$$
If $b(s)$ satisfies the unimodality condition \eqref{eq-UM},
then the above functions are simplified as
$$
b^{*}(u)=\min \left\{b'(0) u,M\right\},
$$
and
$$
b_{\epsilon}(u)=\min\{b(u),d(\zeta_1-\epsilon)\}.
$$
According to the definition, we have
\begin{lemma} \label{le-bstar}
Both $b^*$ and $b_\epsilon$ are continuous on $[0,+\infty)$ and
monotonically increasing on $[0,\zeta_2]$;
$b^*(s)\ge b(s)\ge b_\epsilon(s)$ for all $s\ge0$;
$b^*(\zeta_2)=d(\zeta_2)=M$ and $b^*(s)>d(s)$ for $s\in(0,\zeta_2)$;
$b_\epsilon(\zeta_1-\epsilon)=d(\zeta_1-\epsilon)<d(\zeta_1)$
and $b_\epsilon(s)>d(s)$ for $s\in(0,\zeta_1-\epsilon)$.
\end{lemma}
{\it \bfseries Proof.}
The above statements are obvious and their proofs are omitted for
the sake of simplicity.
$\hfill\Box$

Consider the following two auxiliary
delayed diffusion equations
\begin{equation}\label{eq-zauxupper}
w_{t}(t, x)=D (w^m)_{x x}(t, x)-d(w(t, x))+ b^{*}(w(t-r, y)),
\end{equation}
and
\begin{equation}\label{eq-zauxlower}
w_{t}(t, x)=D (w^m)_{x x}(t, x)-d(w(t, x))+ b_{\epsilon}(w(t-r, y)).
\end{equation}
The wave equations corresponding to
\eqref{eq-zauxupper} and \eqref{eq-zauxlower} are
\begin{equation} \label{eq-zaux1}
c U'(t)-D {U^m}^{\prime \prime}(t)+d(U(t))- b^{*}(U(t-c r))=0,
\end{equation}
and
\begin{equation} \label{eq-zaux2}
c U'(t)-D {U^m}^{\prime \prime}(t)+d(U(t))- b_\epsilon(U(t-c r))=0.
\end{equation}

We note that the characteristic functions of \eqref{eq-zaux1} and \eqref{eq-zaux2}
near $0$ are identical to \eqref{eq-character}, i.e., the
characteristic function of \eqref{eq-tw} near $0$.
However, we will show that the critical wave speed is not determined by
this characteristic function near $0$.

Now we recall the existence of monotone traveling wavefronts for the above
two auxiliary degenerate diffusion equations with time delay.

\begin{lemma}[\cite{JDE18XU}] \label{le-upperNlower}
For any given $m>1$, $D>0$ and $r\ge0$, there exist a constant $\hat c(m,r,b^*,d)$
(depending on $m,r$ and the structure of $b^*(\cdot),d(\cdot)$)
and a constant $\hat c(m,r,b_\epsilon,d)$ (depending on $m,r$ and the
structure of $b_\epsilon(\cdot),d(\cdot)$) such that
\eqref{eq-zaux1} and \eqref{eq-zaux2} admit monotonically increasing wavefronts
$\phi^*(t)$ and $\phi_\epsilon(t)$ for $c_1>\hat c(m,r,b^*,d)$
and $c_2>\hat c(m,r,b_\epsilon,d)$, respectively, with
$0<\phi^*(t)<\zeta_2$, $0<\phi_\epsilon(t)<\zeta_1-\epsilon$,
\begin{equation} \label{eq-leftprofile}
|\phi^*(t)\zeta_2e^{\lambda_{c_1} t}|\le C^*e^{\Lambda_{c_1} t},
\quad |\phi_\epsilon(t)-(\zeta_1-\epsilon)e^{\lambda_{c_2} t}|
\le C_\epsilon e^{\Lambda_{c_2} t},
\quad t<0,
\end{equation}
where $\lambda_1,\lambda_2>0$ are the unique roots of $\chi_0(\lambda)=0$
corresponding to $c_1$ and $c_2$ respectively,
($\chi_0$ is defined in \eqref{eq-character})
and $\Lambda_{c_i}>\lambda_{c_i}$ for $i=1,2$, $C^*,C_\epsilon>0$ are constants.
According to the proof therein, $C_\epsilon$ and $\hat c(m,r,b_\epsilon,d)$
are uniformly bounded to the upwards with respect to $\epsilon\in(0,\epsilon_0)$.
\end{lemma}

{\it \bfseries Proof of Theorem \ref{th-exist}.}
For any given $c>\max\{\hat c(m,r,b^*,d),\hat c(m,r,b_\epsilon,d)\}$,
let $\phi^*(t)$ and $\phi_\epsilon(t)$ be the monotonically increasing wavefronts
of \eqref{eq-zaux1} and \eqref{eq-zaux2}, respectively, corresponding
to the same wave speed $c$.
According to \eqref{eq-leftprofile},
\begin{equation*}
|\phi^*(t)-\zeta_2e^{\lambda t}|\le Ce^{\Lambda t},
\quad |\phi_\epsilon(t)-(\zeta_1-\epsilon)e^{\lambda t}|
\le Ce^{\Lambda t},
\quad t<0,
\end{equation*}
for $\Lambda>\lambda$ and $C>0$ with $\lambda>0$ being
the unique root of $\chi_0(\lambda)=0$ corresponding to $c$.
We may assume that
$$\phi_*(t)\ge \phi_\epsilon(t), \quad \text{~for all ~} t\in\mathbb R.$$
Otherwise, let $t_0<0$ be sufficiently small such that
\begin{equation} \label{eq-zshift}
\phi^*(t)\ge \frac{\zeta_2+\zeta_1}{2}e^{\lambda t}
\ge \phi_\epsilon(t), ~\forall t<t_0,
\end{equation}
and choose $t_1$ such that $\phi^*(t_1)\ge \frac{\zeta_2+\zeta_1}{2}$.
Then we shift $\phi_\epsilon(t)$ to $\phi_\epsilon(t-\max\{t_1-t_0,0\})$.
Therefore,
$$
\phi_\epsilon(t-\max\{t_1-t_0,0\})\le \frac{\zeta_2+\zeta_1}{2}
\le \phi^*(t), \quad \forall t\ge t_1,
$$
and
$$
\phi_\epsilon(t-\max\{t_1-t_0,0\})\le
\phi^*(t-\max\{t_1-t_0,0\})
\le \phi^*(t), \quad \forall t< t_1,
$$
according to \eqref{eq-zshift} and the monotonicity of $\phi^*(t)$.
We replace $\phi_\epsilon(t)$ by $\phi_\epsilon(t-\max\{t_1-t_0,0\})$.

Define
$$
H^{*}[\phi](t)=b^{*}(\phi(t-c r)), \quad t \in \mathbb{R},
$$
and
$$
H_{\epsilon}[\phi](t)=b_\epsilon(\phi(t-c r)), \quad t \in \mathbb{R},
$$
then for any $\phi,\psi\in C(\mathbb R,[0,\zeta_2])$ with $\phi(t)\ge\psi(t),t\in\mathbb R$,
we have
$$
H^{*}[\phi](t) \ge H^{*}[\psi](t)
\quad \text { and } \quad
H_{\epsilon}[\phi](t) \ge H_{\epsilon}[\psi](t) \quad
\text { for all } t \in \mathbb{R},
$$
since $b^*$ and $b_\epsilon$ are monotonically increasing on $[0,\zeta_2]$.
Set
\begin{align} \nonumber
\Gamma_\epsilon:=\Big\{\phi\in C(\mathbb R;\mathbb R);
\phi_\epsilon(t)\le\phi(t)\le\phi^*(t), ~\phi(t)
\text{~is monotonically increasing}
\\ \label{eq-Gamma}
\text{ on~} (-\infty,t_\Gamma],
\text{and~}\phi(t)\ge \phi(t_\Gamma) \text{~for all~} t>t_\Gamma
\Big\},
\end{align}
where $t_\Gamma\in\mathbb R$ is a fixed constant such that
$$
0<\delta_0(\zeta_1-\epsilon_0)\le \phi_\epsilon(t_\Gamma)
\le \phi^*(t_\Gamma)\le \phi^*(t_\Gamma+cr)<\zeta_1
$$
with $\delta_0\in(0,1/2)$ being sufficiently small.
We note that $\phi_\epsilon(t)$ is depending on $\epsilon$,
but the constants in \eqref{eq-leftprofile} can be selected independent
of $\epsilon$ in Lemma \ref{le-upperNlower}, and so is $\delta_0$.
Then we see that $\Gamma_\epsilon$ is nonempty and convex in $\mathscr E$,
where $\mathscr E$ is the linear space $C^b_\text{unif}(\mathbb R)$
endowed with the norm
$$
\|\phi\|_*=\sum_{n=1}^\infty\frac{1}{2^n}\|\phi\|_{L^\infty([-n,n])}.
$$

For any $\psi(t)\in\Gamma_\epsilon$,
we solve the following degenerate equation
\begin{equation} \label{eq-zfixed}
\begin{cases}
\displaystyle
c\phi'(t)-D(\phi^m(t))''+d(\phi(t))=b(\psi(t-cr)), \quad t\in\mathbb R,
\quad \lim_{t\to-\infty}\phi(t)=0, \\
\displaystyle
0<d^{-1}(\liminf_{t\to+\infty}b(\psi(t)))\le
\liminf_{t\to+\infty}\phi(t)
\\
\qquad\qquad
\le\limsup_{t\to+\infty}\phi(t)
\le d^{-1}(\limsup_{t\to+\infty}b(\psi(t)))<+\infty.
\end{cases}
\end{equation}
Denote $$\hat \psi(t):=H[\psi](t):=b(\psi(t-cr)).$$
Since $\psi(t-cr)$ is increasing on $(-\infty,t_\Gamma+cr)$,
$\psi(t)\le \phi^*(t)\le \zeta_1$ for all $t\le t_\Gamma$,
and $b(s)$ is increasing for $s\in[0,\zeta_1]$,
we see that $\hat\psi(t)$ is monotonically increasing on $(-\infty,t_\Gamma]$
and $\hat\psi(t)\ge \hat\psi(t_\Gamma)$ for all $t>t_\Gamma$.
According to Lemma \ref{le-solv},
\eqref{eq-zfixed} admits a solution $\phi(t)$
such that $\phi(t)$ is monotonically increasing on $(-\infty,t_\Gamma]$
and $\phi(t)\ge \phi(t_\Gamma)$ for all $t>t_\Gamma$.
Define $F^*:\Gamma_\epsilon\to C(\mathbb R, [0,\zeta_2])$ by
$F^*(\psi)=\phi$ with $\phi(t)$ being the solution of \eqref{eq-zfixed}
corresponding to $\psi(t)\in\Gamma_\epsilon$.

We need to prove that $F^*(\Gamma_\epsilon)\subset \Gamma_\epsilon$.
For any $\psi(t)\in \Gamma_\epsilon$, we have
$\phi_\epsilon(t)\le\psi(t)\le\phi^*(t)$, then
$$
H[\psi](t-cr)\le H^*[\psi](t-cr) \le H^*[\phi^*](t-cr),
$$
and
\begin{equation*}
\begin{cases}
\displaystyle
c\phi^*{}'(t)-D(\phi^{*m}(t))''+d(\phi^*(t))
\ge c\phi'(t)-D(\phi^m(t))''+d(\phi(t)), \quad t\in\mathbb R, \\
\displaystyle
\liminf_{t\to-\infty}(\phi^*(t)-\phi(t))
=\liminf_{t\to-\infty}\phi^*(t)-\liminf_{t\to-\infty}\phi(t)=0, \\
\displaystyle
\limsup_{t\to+\infty}(\phi^*(t)-\phi(t))
\ge\lim_{t\to+\infty}\phi^*(t)-\limsup_{t\to-\infty}\phi(t) \\
\displaystyle
\qquad\qquad\ge\zeta_2-d^{-1}(\limsup_{t\to+\infty}b(\psi(t)))
\ge\zeta_2-d^{-1}(\limsup_{t\to+\infty}b^*(\phi^*(t)))=0,
\end{cases}
\end{equation*}
since $\phi(t)$ and $\phi^*(t)$ are solutions of \eqref{eq-zfixed}
and \eqref{eq-zaux1}.
Applying the comparison principle Lemma \ref{le-comparison},
we find $\phi(t)\le \phi^*(t)$ for all $t\in\mathbb R$.
In a similar way, the property $\phi(t)\ge \phi_\epsilon(t)$
follows from the comparison principle Lemma \ref{le-comparison}
and the inequality
$$
H[\psi](t-cr)\ge H_\epsilon[\psi](t-cr) \ge H_\epsilon[\phi_\epsilon](t-cr).
$$

From the proof of Lemma \ref{le-solv}, we see that the solutions $\phi(t)$
of \eqref{eq-zfixed} are uniformly bounded in $C^\alpha([-n,n])$
with some $\alpha\in(0,1/(2m))$,
$\phi^m(t)$ are uniformly bounded in $W^{1,2}([-n,n])$ for
any compact interval $[-n,n]$,
and $\phi(t)$ are uniformly bounded in $L^\infty(\mathbb R)$.
According to the definition of the function space $\mathscr E$,
$F^*(\Gamma_\epsilon)$ is compact in $\mathscr E$.
By the Schauder's fixed point theorem,
it follows that $F^*$ has a fixed point
$U$ in $\Gamma_\epsilon\subset \mathscr E$,
which satisfies
$$
c U'(t)-D {U^m}^{\prime \prime}(t)+d(U(t))- b(U(t-c r))=0,
$$
and
\begin{equation}\label{eq-thexepsi}
\phi_\epsilon(t)\le U(t)\le \phi^*(t) \quad \text{for all}\quad  t\in\mathbb R.
\end{equation}
Moreover, $U(-\infty)=0$ and
$$
\zeta_1-\epsilon \le \liminf _{t \rightarrow+\infty} U(t)
\le \limsup _{t \rightarrow+\infty} U(t) \le \zeta_2.
$$
Since $U(t)$ is independent of $\epsilon$, taking the limit as
$\epsilon\to 0^+$, we have
$$
\zeta_1\le \liminf _{t \rightarrow+\infty} U(t)
\le \limsup _{t \rightarrow+\infty} U(t) \le \zeta_2.
$$
The proof is completed.
$\hfill\Box$

\section{Nonexistence of traveling wave solutions}

This section is devoted to the proof of Theorem \ref{th-non-exist}.
The proof is based on the phase transform approach
similar to the proof of Lemma 3.11 in \cite{JDE18XU}
with some modification suitable for large time delay
and non-monotone birth rate functions.
We note that this method is incapable of showing the existence of
traveling waves with time delay in general
since the trajectories with time delay may intersect with each other.
However, it can be a blueprint to draw a contradiction for
proving the nonexistence.

{\it \bfseries Proof of Theorem \ref{th-non-exist}.}
For any given $m>1$, $D>0$ and $r\ge0$, we need to find a constant
$\dot c(m,r,b,d)>0$, such that, \eqref{eq-tw} admits no wave solution $\phi(t)$
(semi-wavefronts or wavefronts, sharp or smooth)
for any $c<\dot c(m,r,b,d)$.
The nonexistence result is valid
for a typical Nicholson's birth rate function and death rate function
without time delay in \cite{Huang-Jin-Mei-Yin}.
We can verify that it is also true for the general type of $b$ and $d$
without time delay.
Here we only prove the case with time delay $r>0$.

We prove by contradiction and assume that $\phi_c$ is a wave solution
corresponding to the speed $c$.
Recall that $\zeta_1$ and $\zeta_2$ are the constants in \eqref{eq-zeta}.
Since $b'(0)>d'(0)$, let $(0,\zeta_3)$ be the maximal interval
such that
\begin{equation} \label{eq-zpsi0}
\psi_0(\phi):=\frac{Dm\phi^{m-1}(b(\phi)-d(\phi))}{c}
\end{equation}
is increasing with respect to $\phi$
and denote $\zeta_0=\min\{\zeta_1,\zeta_3\}$.
It should be noted that $\zeta_0$ is independent of $c$ and $r$.
Let $I_0:=(-\infty,t_0)$ be the maximal interval of the leading edge of $\phi_c$
and let $I_1=(t_1,t_2)$ be the maximal subinterval of $I_0$ such that $\phi_c$
is positive, monotonically increasing and $\phi_c(t)<\zeta_0$.
That is, there exists a unique $\hat t_0<t_0$ such that $\phi_c(\hat t_0)=\zeta_0$
and we take $t_2=\hat t_0$.
If $\phi_c$ is of smooth type, then $t_1=-\infty$,
while if $\phi_c$ is sharp, we take $t_1=0$ instead.
Within $I_1$, $\phi_c(t)$ is monotonically increasing and
$b(\phi_c)$ is monotonically increasing with respect to $\phi_c$
as $\phi_c\le \zeta_0\le\zeta_1$.

Now we introduce the phase transform approach, see for example
\cite{Huang-Jin-Mei-Yin,JDE18XU}.
Let
$$\psi_c(t)=D(\phi_c^m(t))'.$$
Since $\phi_c(t)$ is positive and monotonically increasing in $I_1$,
we have the following singular phase plane
where $(\phi_c(t),\psi_c(t))$ corresponds to a trajectory
\begin{equation} \label{eq-nonex-pp}
\begin{cases}
\displaystyle
\phi'(t)=\frac{\psi(t)}{Dm\phi^{m-1}(t)}=:\Phi,\\
\displaystyle
\psi'(t)=\frac{c\psi(t)}{Dm\phi^{m-1}(t)}+d(\phi(t))-b(\phi_{cr}(t)):=\Psi,
\end{cases}
\end{equation}
with $\phi_{cr}(t)=\phi(t-cr)$.
We note that $\psi_c(t)\ge0$ and the zero points of $\psi_c(t)$ is isolated
since $\phi_c(t)$ is a given wave solution.
According to the choice of $I_1$,
we can regard $t\in I_1$ as
a inverse function of $\phi_c$ and denote
$\tilde\psi_c(\phi_c)=\psi_c(t(\phi_c))\ge0$.
We redefine $\phi_{cr}(t)$ as a functional of $\phi_c$ and $\tilde\psi_c$ as follows
\begin{equation} \label{eq-zdelaypp}
\phi_{cr}=\inf_{\theta\in[0,\phi_c]}\Big\{\int_\theta^{\phi_c}
\frac{Dms^{m-1}}{\tilde\psi_c(s)}ds\le cr\Big\}.
\end{equation}
Consider the following nonlocal problem
\begin{equation} \label{eq-zpp}
\begin{cases}
\displaystyle
\frac{d\psi}{d\phi}=c-\frac{Dm\phi^{m-1}(b(\phi_{cr})-d(\phi))}{\psi}
=\frac{\Psi}{\Phi}, \\
\psi(0)=0, \quad \psi(\zeta_0)=Dm\zeta_0^{m-1}\phi_c'(t_2)>0,
\qquad \phi\in(0,\zeta_0).
\end{cases}
\end{equation}
Here, nonlocal means that $\phi_{cr}$ is a functional of $\phi$ and $\psi(\phi)$,
which is caused by the time delay.

We draw a contradiction to the existence of solutions to \eqref{eq-zpp}
when $c$ is sufficiently small with the help of the phase plane \eqref{eq-nonex-pp}.
The curve $\Gamma_c$ corresponding to $\psi_0(\phi)$
defined in \eqref{eq-zpsi0} divides
$(0,\zeta_0)\times(0,+\infty)$ into two parts,
$E_1:=\{(\phi,\psi);\phi\in(0,\zeta_0),0<\psi<\psi_0(\phi)\}$
and $E_2:=((0,\zeta_0)\times(0,+\infty))\backslash E_1$.
For any $(\phi,\psi)\in\Gamma_c$, we have
$$
\frac{\Psi}{\Phi}=c-\frac{Dm\phi^{m-1}(b(\phi_{cr})-d(\phi))}{\psi}
>c-\frac{Dm\phi^{m-1}(b(\phi)-d(\phi))}{\psi}=0.
$$
We can check that $\Psi/\Phi>0$ for any $(\phi,\psi)\in E_2$.
Let $\psi_1(\phi)$ be the solution of
$$
\begin{cases}
\displaystyle
\frac{d\psi}{d\phi}=c+\frac{Dm\phi^{m-1}d(\phi)}{\psi}, \\
\psi(0)=0, \psi(\phi)>0, \phi\in(0,\zeta_0).
\end{cases}
$$
Asymptotic analysis shows that there exists a constant $C_1>0$
depending on the upper bound of $c$ (independent of $c$ if $c$ is small) such that
$$\phi_1(\phi)\le C_1\phi, \quad \phi\in(0,\zeta_0).$$
The comparison principle of \eqref{eq-zpp} shows that
\begin{equation} \label{eq-ztildepsi}
\tilde\psi_c(\phi)\le \phi_1(\phi)\le C_1\phi, \quad \phi\in(0,\zeta_0).
\end{equation}
Let $\epsilon\in(0,\zeta_0)$ be a constant such that
\begin{equation} \label{eq-zeps}
\int_0^\epsilon \phi^{m-1}d(\phi)d\phi<
\frac{1}{4}\int_\epsilon^{\zeta_0}\phi^{m-1}(b(\phi)-d(\phi))d\phi,
\end{equation}
and
$$
\delta:=\inf_{\phi\in(\epsilon,\zeta_0)}(b(\phi)-d(\phi))>0.
$$
We note that $\epsilon$ and $\delta$ only depend on the structure of $b$ and $d$.
We assert that for any given $r>0$, there exists a $c_1>0$ such that
$b(\phi_{cr})-d(\phi)>0$ for all $\phi\in(\epsilon,\zeta_0)$ if $c\le c_1$.
In fact, according to \eqref{eq-ztildepsi},
$$
c_1r\ge
cr=\int_{\phi_{cr}}^\phi \frac{Dms^{m-1}}{\tilde\psi_c(s)}ds
\ge \int_{\phi_{cr}}^\phi \frac{Dms^{m-1}}{C_1s}ds
\ge \frac{Dm}{C_1(m-1)}(\phi^{m-1}-\phi_{cr}^{m-1}),
$$
and then using the uniform continuity of the function $f(s):=s^{1/(m-1)}$
on the interval $[\epsilon/2,\zeta_0]$
with the continuity modulus function being denoted by $\omega(\cdot)$, we have
\begin{align*}
0<b(\phi)-b(\phi_{cr})&=b'(\theta)(\phi-\phi_{cr})
\le \sup_{s\in(0,\zeta_0)}b'(s)\cdot(\phi-\phi_{cr})
\\
&\le \sup_{s\in(0,\zeta_0)}b'(s)\cdot\omega(|\phi^{m-1}-\phi_{cr}^{m-1}|)
\\
&\le \sup_{s\in(0,\zeta_0)}b'(s)\cdot\omega(\frac{C_1(m-1)c_1r}{Dm})
\le \frac{\delta}{2}
\end{align*}
for some $\theta\in(\phi_{cr},\phi)$,
provided that $c_1r$ is sufficiently small
such that $c_1r=\mu_0:=\mu_0(m,b(\cdot),d(\cdot))>0$
(it suffices that $c_1$ is sufficiently small as $r$ is given).
Here we note that $\mu_0(m,b(\cdot),d(\cdot))$ is a constant
depending on $m$, $\epsilon$, $\zeta_0$, $\delta$, $\sup_{s\in(0,\zeta_0)}b'(s)$,
which are all dependent on $m$ and the structure of $b(\cdot)$ and $d(\cdot)$.
The dependence of $\mu_0(m,b(\cdot),d(\cdot))$ on $b(\cdot)$ is basically on
the structure of $b(\cdot)$ within $(0,\zeta_0)$ and
$\zeta_0\le \zeta_1$ with $\zeta_1$ depending on the whole structure of
$b(\cdot)$ on $(0,\zeta_2)$.
Therefore,
\begin{align} \nonumber
b(\phi_{cr})-d(\phi)&=(b(\phi)-d(\phi))-(b(\phi)-b(\phi_{cr}))
\\ \label{eq-zbphicr}
&\ge (b(\phi)-d(\phi))-\frac{\delta}{2}
\ge \frac{b(\phi)-d(\phi)}{2},
\quad \phi\in(\epsilon,\zeta_0).
\end{align}
The first integral of \eqref{eq-zpp} over $(0,\zeta_0)$ shows that
\begin{align*}
c\int_0^{\zeta_0}\tilde\psi_c(\phi)d\phi
&=\frac{1}{2}\tilde\psi_c^2(\phi)\Big|_0^{\zeta_0}+
\int_0^{\zeta_0} Dm\phi^{m-1}(b(\phi_{cr})-d(\phi))d\phi\\
&\ge \int_0^\epsilon Dm\phi^{m-1}(b(\phi_{cr})-d(\phi))d\phi
+\int_\epsilon^{\zeta_0} Dm\phi^{m-1}(b(\phi_{cr})-d(\phi))d\phi\\
&\ge -\int_0^\epsilon Dm\phi^{m-1}d(\phi)d\phi
+\int_\epsilon^{\zeta_0} Dm\phi^{m-1}(b(\phi_{cr})-d(\phi))d\phi\\
&\ge \Big(-\frac{1}{4}+\frac{1}{2}\Big)\int_\epsilon^{\zeta_0} Dm\phi^{m-1}(b(\phi)-d(\phi))d\phi,
\end{align*}
where we have used \eqref{eq-zeps} and \eqref{eq-zbphicr}.
On the other hand, we have
$$
c\int_0^{\zeta_0}\tilde\psi_c(\phi)d\phi
\le c\int_0^{\zeta_0}C_1\phi d\phi
\le c\frac{C_1}{2}\zeta_0^2.
$$
Now we arrive at a contradiction if we have chosen
$c\le \dot c$ with
$$\dot c=\min\{c_1,c_2\}=\min\{\frac{\mu_0(m,b(\cdot),d(\cdot))}{r},
c_2(m,b(\cdot),d(\cdot))\}$$
such that
$$
c_2\frac{C_1}{2}\zeta_0^2<\frac{1}{4}\int_\epsilon^{\zeta_0} Dm\phi^{m-1}(b(\phi)-d(\phi))d\phi.
$$
The proof is completed.
$\hfill\Box$

\section{Existence of sharp waves}

In this section,
we develop a new delayed iteration approach
based on an observation of the delicate structure of time delay and sharp edge
to solve the delayed degenerate equation.
As far as we know, this is the first framework
of showing the existence of sharp traveling wave solution for the degenerate
diffusion equation with large time delay.
A sharp wave solution $\phi(t)$ is a special solution such that
$\phi(t)\equiv0$ for $t\le0$ and $\phi(t)>0$ for $t>0$.
The existence of sharp wave solution for the case without time delay
and with Nicholson's birth rate function $b(u)=pue^{-au}$ and death rate function
$d(u)=\delta u$ for some constants $p,a,\delta$ is proved in \cite{Huang-Jin-Mei-Yin}.
It is also valid for the general birth rate and death rate functions without time delay
and here we only focus on the case with time delay.

For any given $m>1$, $D>0$ and $r>0$,
we solve \eqref{eq-tw} step by step.
First, noticing that the sharp wave solution $\phi(t)=0$ for $t\le0$
and then $\phi(t-cr)=0$ for $t\in[0,cr)$, \eqref{eq-tw} is locally reduced to
\begin{equation} \label{eq-semi-1}
\begin{cases}
c\phi'(t)=D(\phi^m(t))''-d(\phi(t)), \\
\phi(0)=0, \quad (\phi^m)'(0)=0, \quad t\in(0,cr),
\end{cases}
\end{equation}
whose solutions are not unique and we choose the maximal one
such that $\phi(t)>0$ for $t\in(0,cr)$ as shown in the following lemma.
Here, $(\phi^m)'(0)=0$ is necessary and sufficient condition
such that the zero extension of $\phi(t)$ to the left
satisfies \eqref{eq-tw} locally near $0$ in the sense of distributions.

\begin{lemma} \label{le-semi-1}
For any $c>0$, the degenerate ODE \eqref{eq-semi-1} admits a maximal solution
$\phi_c^1(t)$ on $(0,cr)$ such that $\phi_c^1(t)>0$ on $(0,cr)$ and
$$
\phi_c^1(t)=\Big(\frac{(m-1)c}{Dm}t\Big)^\frac{1}{m-1}+o(t^\frac{1}{m-1}),
\quad t\to0^+.
$$
\end{lemma}
{\it\bfseries Proof.}
Clearly, $\phi_0(t)\equiv0$ is a solution of \eqref{eq-semi-1}.
But we are looking for the solution such that $\phi_c^1(t)>0$ on $(0,cr)$.
Consider the generalized phase plane related to \eqref{eq-semi-1} and
define $\psi_c^1(t)=D[(\phi_c^1(t))^m]'$, then $(\phi_c^1(t),\psi_c^1(t))$
solve the following singular ODE system on $(0,cr)$
\begin{equation} \label{eq-semi-pp-1}
\begin{cases}
\displaystyle
\phi'(t)=\frac{\psi(t)}{Dm\phi^{m-1}(t)},\\
\displaystyle
\psi'(t)=\frac{c\psi(t)}{Dm\phi^{m-1}(t)}+d(\phi(t)).
\end{cases}
\end{equation}
We solve \eqref{eq-semi-pp-1} with the condition
$(\phi_{c,\epsilon}^1(0),\psi_{c,\epsilon}^1(0))=(0,\epsilon)$ with $\epsilon>0$,
whose existence, continuous dependence and suitable regularity
follow from the phase plane analysis.
Let $\epsilon$ tends to zero and $(\phi_c^1(t),\psi_c^1(t))$ be
the limiting function.
Then $\phi_c^1(t)$ is the maximal solution of \eqref{eq-semi-1}
and $\phi_c^1(t)>0$ on $(0,cr)$.
Asymptotic analysis shows that
$$
\psi_c^1(t)=Dm(\phi_c^1(t))^{m-1}\phi_c^1{}'(t)
=c\phi_c^1(t)+o(\phi_c^1(t)), \quad t\to0^+,
$$
which means that
$$
\phi_c^1(t)=\left(\frac{(m-1)c}{Dm}t+o(t)\right)^\frac{1}{m-1},
\quad t\to0^+.
$$
$\hfill\Box$

Second, let $\phi_c^2(t)$ be the solution of the following initial value ODE problem
\begin{equation} \label{eq-semi-2}
\begin{cases}
c\phi'(t)=D(\phi^m(t))''-d(\phi(t))+b(\phi_c^1(t-cr)), \\
\phi(r)=\phi_c^1(r), \quad
\phi'(r)=(\phi_c^1)'(r), \qquad t\in(cr,2cr).
\end{cases}
\end{equation}
Define $\psi_c^2(t)=D[(\phi_c^2(t))^m]'$, then $(\phi_c^2(t),\psi_c^2(t))$
solve the following system on $(cr,2cr)$
\begin{equation} \label{eq-semi-pp-2}
\begin{cases}
\displaystyle
\phi'(t)=\frac{\psi(t)}{Dm\phi^{m-1}(t)},\\
\displaystyle
\psi'(t)=\frac{c\psi(t)}{Dm\phi^{m-1}(t)}+d(\phi(t))-b(\phi_c^1(t-cr)).
\end{cases}
\end{equation}
The above steps can be continued unless $\phi_c^k(t)$ blows up
or decays to zero in finite time for some $k\in \mathbb N^+$.
Let $\phi_c(t)$ be the connecting function of those functions on each step, i.e.,
\begin{equation} \label{eq-semi}
\phi_c(t)=
\begin{cases}
\phi_c^1(t), \quad &t\in[0,cr),\\
\phi_c^2(t), \quad &t\in[cr,2cr),\\
\dots\\
\phi_c^k(t), \quad &t\in[(k-1)cr,kcr),\\
\dots
\end{cases}
\end{equation}
for some finite steps such that $\phi_c(t)$ blows up or decays to zero,
or for infinite steps such that $\phi_c(t)$ is defined on $(0,+\infty)$
and zero extended to $(-\infty,0)$ for convenience.

\begin{lemma} \label{le-semi-decay}
For any given $m$, $D$ and $r>0$, there exists a constant $\underline c>0$ such that
if $c\le \underline c$, then
$\phi_c(t)$ decays to zero in finite time.
\end{lemma}
{\it\bfseries Proof.}
On the existence interval of $\phi_c(t)$, the pair
$(\phi_c(t),\psi_c(t))$ with $\psi_c(t):=D[(\phi_c(t))^m]'$
is a trajectory in the phase plane
\begin{equation} \label{eq-zsemipp}
\begin{cases}
\displaystyle
\phi'(t)=\frac{\psi(t)}{Dm\phi^{m-1}(t)},\\
\displaystyle
\psi'(t)=\frac{c\psi(t)}{Dm\phi^{m-1}(t)}+d(\phi(t))-b(\phi(t-cr)).
\end{cases}
\end{equation}
The proof of $\phi_c(t)$ decays to zero in finite time
is similar to the proof of nonexistence of
semi-wavefront with monotonically increasing leading edge when $c$ is sufficiently small,
i.e., the proof of Theorem \ref{th-non-exist}.
Here we omit the proof.
$\hfill\Box$

\begin{lemma} \label{le-semi-blow}
For any given $m$, $D$ and $r>0$, there exists a constant $\overline c>0$ such that
if $c\ge \overline c$, then
$\phi_c(t)$ grows up to $+\infty$ as $t$ tends to $+\infty$.
\end{lemma}
{\it\bfseries Proof.}
On the existence interval of $\phi_c(t)$, the pair
$(\phi_c(t),\psi_c(t))$ defined in the proof of Lemma \ref{le-semi-decay}
is a trajectory in the phase plane \eqref{eq-zsemipp}.
Now, we utilize the phase plane analysis to show that when $c$ is large enough,
then $\phi_c(t)$ grows up to the positive infinity as $t$ increases.
For $t\in(0,cr)$, $\phi_c$ is strictly monotonically increasing
according to \eqref{eq-semi-pp-1}.
Let $(0,\zeta)$ be the maximal interval such that $\phi_c$ is
strictly monotonically increasing and within this interval, we have
${d\psi_c}/{d\phi_c}$ satisfies
\begin{equation} \label{eq-zsemi-blow}
\begin{cases}
\displaystyle
\frac{d\psi}{d\phi}=c-\frac{Dm\phi^{m-1}(b(\phi_{cr})-d(\phi))}{\psi}
=:\frac{\Psi}{\Phi}, \\
\psi(0)=0, \quad \psi(\phi)>0,
\qquad \phi\in(0,\zeta),
\end{cases}
\end{equation}
as in the proof of Theorem \ref{th-non-exist},
where $\phi_{cr}$ is the functional of $\phi_c$
and $\psi_c$ defined in \eqref{eq-zdelaypp}
(we regard $\psi_c$ as a function of $\phi_c$
since $\phi_c$ is strictly increasing).
Consider the following auxiliary problem
\begin{equation} \label{eq-zsemi-blow-aux}
\begin{cases}
\displaystyle
\frac{d\psi}{d\phi}=c-\frac{Dm\phi^{m-1}(\tilde b(\phi)-d(\phi))}{\psi}, \\
\psi(0)=0, \quad \psi(\phi)>0,
\qquad \phi\in(0,\zeta),
\end{cases}
\end{equation}
where $\tilde b(s)=\sup_{\theta\in(0,s)}b(\theta)$
is the quasi-monotone modification of $b(s)$
and the solution of \eqref{eq-zsemi-blow-aux} is denoted by
$\underline \psi{}_c(\phi)$.
Therefore, as $\phi_c(t)$ is strictly increasing (equivalently, $\psi_c(t)>0$) we have
$$
b(\phi_{cr})\le \tilde b(\phi_{cr}) \le \tilde b(\phi),
$$
and the comparison between \eqref{eq-zsemi-blow} and \eqref{eq-zsemi-blow-aux} shows that
\begin{equation} \label{eq-zsemi-com}
\psi_c(\phi)\ge \underline\psi{}_c(\phi), \qquad \phi\in(0,\zeta).
\end{equation}
The phase plane analysis to \eqref{eq-zsemi-blow-aux} without time delay
shows that there exists a $\overline c>0$ such that
if $c\ge \overline c$, then
$\underline\psi{}_c(\phi)$ is positive for all $\phi\in(0,+\infty)$,
which means according to \eqref{eq-zsemi-com} that
$\psi_c(\phi)>0$ for all $\phi\in(0,+\infty)$,
$\phi_c(t)$ is always increasing for $t\in(0,+\infty)$.
It follows that in fact $\zeta=+\infty$ and $\phi_c(t)$ grows up to $+\infty$
as $t$ tends to $+\infty$.
$\hfill\Box$

We also need the following continuous dependent property of $\phi_c(t)$ on $c$
proved in \cite{Our-Forthcoming}.

\begin{lemma}[\cite{Our-Forthcoming}] \label{le-continu}
For any given $m$, $D$ and $r>0$,
the solution $\phi_c(t)$ is locally continuously dependent on $c$.
That is, for any $c>0$ and any given $T>0$ and $\varepsilon>0$,
there exists a $\delta>0$ such that for any $|c_1-c|<\delta$ and $c_1>0$ we have
$$
|\phi_{c_1}(t)-\phi_c(t)|<\varepsilon,
\quad \forall t\in(0,T_1-\varepsilon),
$$
where $T_1=\min\{T,T_c\}$
with $T_c$ being the existence interval of $\phi_c(t)$.
\end{lemma}

Now, we are able to prove the existence of sharp traveling waves.

{\it \bfseries Proof of Theorem \ref{th-semifinite}.}
Let $(0,T_1)$ and $(0,T_2)$ be the maximal interval such that
$\phi_{\underline c}(t)$ remains positive before decaying to zero
and $\phi_{\overline c}(t)<\zeta_2$, respectively,
where $\underline c$ and $\overline c$ are constants in Lemma \ref{le-semi-decay}
and Lemma \ref{le-semi-blow}.
For any $T>\max\{T_1,T_2\}$, $\phi_c(T)\ge\zeta_2$ for some $c\ge \overline c$
and $\phi_c(T)\le0$ for some $c\le \underline c$.
The continuous dependence of $\phi_c(t)$ with respect to $c$ on
the compact interval $[0,T]$ (Lemma \ref{le-continu}) implies that
there exists a $c_T\in[\underline c,\overline c]$ such that
$\phi_{c_T}(T)=\kappa$.
Since the closed interval $[\underline c,\overline c]$ is compact,
there exists a subsequence of $\{c_T\}$, i.e., $\{c_{T_i}\}_{i=1}^\infty$,
and a $c_0\in [\underline c,\overline c]$,
such that $\lim_{i\to\infty}c_{T_i}=c_0$.
Meanwhile, $\phi_{c_0}(t)$ exists on the whole $(0,+\infty)$,
whose zero extension to the left is a sharp wave solution.
The uniform permanence property Lemma \ref{le-permanence}
and the asymptotic expansion Lemma \ref{le-semi-1}
indicate that the sharp wave solutions $\phi_{c_0}(t)$ satisfies
$$0<\zeta_1\le\liminf_{t\to+\infty}\phi_{c_0}(t)
\le\limsup_{t\to+\infty}\phi_{c_0}(t)\le\zeta_2,$$
and
$$
|\phi_{c_0}(t)-C_1t_+^\lambda|\le C_2t_+^\Lambda,
\quad \text{~for any~} t\in(0,1),
$$
where $t_+=\max\{t,0\}$, $\lambda=1/(m-1)$
and $\Lambda>\lambda$, $C_1,C_2>0$ are constants.
$\hfill\Box$

\begin{remark}
The time delay together with the non-monotone structure
of birth rate function $b(u)$ causes us
essential difficulty in proving the monotonic dependence of $\phi_c(t)$
with respect to $c$.
Actually, the possible existence of non-monotone semi-wavefront
suggests that the monotonic dependence may be violated in general.
Without this monotonic dependence, the uniqueness of the
wave speed for wave solutions of sharp type remains open.
\end{remark}

{\it\bfseries Proof of Theorem \ref{th-sharp}.}
The asymptotic behavior near $0$ in Lemma \ref{le-semi-1}
completes the proof.
$\hfill\Box$

\section{Traveling wave solutions with oscillations}

In this section, we follow the main line of \cite{TrofimchukJDE}
to show the oscillating of the wave solutions.
The monotonicity or oscillating, convergence or non-decaying oscillation,
are the basic features of the asymptotic behavior for the wave solutions near
the positive equilibrium $\kappa$.
We note that the nonlinear diffusion equation \eqref{eq-tw}
does not degenerate near $\kappa$
and its linearization near $\kappa$ is of the same type
as the linear diffusion case.
Those observations made us enable to apply the method in \cite{TrofimchukJDE}
and \cite{Gomez} to our nonlinear diffusion case.

Here we recall the concept of slowly oscillating solutions of \eqref{eq-tw},
see for example \cite{TrofimchukJDE}.

\begin{definition}
Let $\psi:[\theta,+\infty) \rightarrow \mathbb{R}$ be a continuous function
for some $\theta\in\mathbb R$.
We say that $\psi$ is oscillatory if
there exist sequences $\{t_{n}\}_{n \ge 1}$ and $\{t_{n}'\}_{n \ge 1}$
such that $t_{n}, t_{n}' \rightarrow+\infty$ and $\psi(t_{n})<0<\psi(t_{n}'), n \ge 1$.
\end{definition}

\begin{definition}[\cite{TrofimchukJDE}]
Set $\mathbb{K}=[-r,0] \cup\{1\}$.
For any $v \in C(\mathbb{K}\backslash\{0\})$
we define the number of sign changes by
$$
\operatorname{sc}(v)=\sup \left\{k \ge 1:
\text { there are } t_{0}<\cdots<t_{k} \text { such that } v(t_{i-1}) v(t_{i})<0
\text { for } i \ge1\right\}.
$$
We set $\operatorname{sc}(v)=0$ if $v(s) \ge 0$ or $v(s) \le 0$ for $s \in \mathbb{K}.$ If $\varphi :[a-r,+\infty) \rightarrow \mathbb{R}$ is a solution of
\eqref{eq-tw}, we set
$\left(\overline{\varphi}_{t}\right)(s)=\varphi(t+s)-\kappa$ if $s \in[-r, 0],$ and $\left(\overline{\varphi}_{t}\right)(1)=\varphi'(t).$
We will say that $\varphi(t)$ is slowly oscillating about $\kappa$ if $\varphi(t)-\kappa$ is oscillatory and for each $t\ge a$, we have either $\operatorname{sc}\left(\overline{\varphi}_{t}\right)=1$ or $\operatorname{sc}\left(\overline{\varphi}_{t}\right)=2$.
\end{definition}

The characteristic function near $\kappa$ plays an essential role in the
investigation of the monotonicity near $\kappa$.
Since the linearization of \eqref{eq-tw} near $\kappa$ is of the same type of
the linear diffusion case, we have the following results
as Lemma 1.1 in \cite{Gomez}.

\begin{lemma}[\cite{Gomez}] \label{le-characteristick}
For $b'(\kappa)<0$,  there exists an extended real number
$c_{\kappa}=c_{\kappa}(m,r,b'(\kappa),d'(\kappa)) \in(0,+\infty]$
such that the characteristic equation
$\chi_\kappa(\lambda)$ defined in \eqref{eq-character-k}
has three real roots $\lambda_{1} \le \lambda_{2}<0<\lambda_{3}$
if and only if $c \le c_{\kappa}.$
If $c_{\kappa}$ is finite and $c=c_{\kappa},$ then $\chi_\kappa(\lambda)$
has a double root $\lambda_{1} =\lambda_{2}<0$,
while for $c>c_{\kappa}$ there does not exist any negative
root to \eqref{eq-character-k}.
Moreover, if $\lambda_{j} \in \mathbb{C}$ is a complex root of \eqref{eq-character-k}
for $c \in(0, c_{\kappa}],$ then $\Re \lambda_{j}<\lambda_{2}$.
Furthermore, $c_\kappa(m,0,b'(\kappa),d'(\kappa))=+\infty$
and $c_\kappa$ is strictly decreasing in its domain,
$$
c_\kappa(m,r,b'(\kappa),d'(\kappa))=
\frac{\mu_\kappa(m,b'(\kappa),d'(\kappa))+o(1)}{r},
\quad r\to+\infty,
$$
where $\mu_\kappa(m,b'(\kappa),d'(\kappa)):=
\sqrt{\frac{2Dm\kappa^{m-1}\omega_\kappa}{b'(\kappa)}}e^{\frac{\omega_\kappa}{2}}$,
and $\omega_\kappa<0$ is the unique negative root of
$2d'(\kappa)=b'(\kappa)e^{-\omega_\kappa}(2+\omega_\kappa)$.
\end{lemma}

We also need the following auxiliary result,
which is Corollary 24 in \cite{TrofimchukJDE}.

\begin{lemma}[\cite{TrofimchukJDE}]\label{le-preosc}
Assume that $f : \mathbb{R}_{+} \rightarrow \mathbb{R}_{+}, f(+\infty)=0$,
does not decay superexponentially.
Then for every $\rho>0$, there exist a sequence
$t_{j} \rightarrow+\infty$ and a real $\delta>1$ such that
$f(t_{j})=\max _{s \ge t_{j}} f(s)$ and
$\max_{s \in[t_{j}-\rho, t_{j}]} f(s) \le \delta f(t_{j})$.
\end{lemma}

Now we prove that the semi-wavefronts are oscillating if $c>c_\kappa$
in a similar method as Lemma 25 in \cite{TrofimchukJDE}
and Lemma 4.6 in \cite{Gomez}.

\begin{lemma}\label{le-notmonotone}
Assume that $b'(\kappa)<0$ and $c>c_{\kappa}$ as in Lemma \ref{le-characteristick},
then \eqref{eq-tw} does not have any
eventually monotone semi-wavefront.
\end{lemma}
{\it\bfseries Proof.}
The proof is similar to the one of Lemma 25 in \cite{TrofimchukJDE}.
Here we provide a sketch of proof using slightly different arguments
suitable for nonlinear diffusion.
Lemma \ref{le-characteristick} implies that
the characteristic function $\chi_\kappa(\lambda)$ around $\kappa$
does not have any negative zeros.
Arguing by contradiction, suppose that, there exists an
eventually monotone travelling wave front.

Set $w(t)=\phi(t)-\kappa,$ then $w(t)$ is either decreasing and strictly positive or increasing and
strictly negative on some interval $[T,+\infty)$ and satisfies
\begin{equation}\label{eq-w1}
Dm(\phi(t)^{m-1}w'(t))'-cw'(t)=p(t)w(t)+k(t)w(t-h),
\end{equation}
where $h=cr$ and
$$
k(t):=-\frac{b(\phi(t-h))-b(\kappa)}{\phi(t-h)-\kappa},
\quad p(t):=\frac{d(\phi(t))-d(\kappa)}{\phi(t)-\kappa}.
$$
Since $\phi(+\infty)=\kappa,$
$0<k(t)<-2 b'(\kappa),$
and $0<p(t)<2d'(\kappa)$ for all sufficiently large $t$.
We will show that for $c>c_{\kappa}$, $w(t)$ will oscillate about zero.
As a consequence of Lemma 3.1.1 from \cite{Hupke},
we can conclude that $w(t)$ cannot convey superexponentially to $0$.
This fact and Lemma \ref{le-preosc} imply the existence of a sequence
$t_{j} \rightarrow+\infty$ and a real number $\delta>0$ such that
$|w(t_{j})|=\max_{s \ge t_{j}}|w(s)|$ and
$\max_{s \in[t_{j}-3 h, t_{j}]}|w(s)| \le \delta|w(t_{j})|$ for every $j.$
Without loss of generality we assume that $w'(t_{n}) \le 0$
and $0<w(t) \le w(t_{n})$ for all $t \ge t_{n}$.
Additionally, we can find a sequence $\{s_{j}\}$ with
$\lim (s_{j}-t_{j})=+\infty$ such that $|w'(s_{j})| \le w(t_{j}).$
Now, since $w(t)$ satisfies \eqref{eq-w1}, we conclude that every
$y_{j}(t)=w(t+t_{j})/w(t_{j})>0$ is a solution of
$$
Dm(\phi^{m-1}(t+t_j)y')'-cy'-p(t+t_j)y-k(t+t_j)y(t-h)=0,
$$
It is clear that $\lim_{j \rightarrow+\infty} k(t+t_{j})=-b'(\kappa)$,
$\lim_{j \rightarrow+\infty} p(t+t_{j})=d'(\kappa)$,
and $\lim_{j\to\infty}\phi(t+t_j)=\kappa$ uniformly on $\mathbb{R}_{+}$
and also that $0<y_{j}(t) \le \delta$ for all
$t \ge-3 h, j=1,2,3, \ldots$.

We need to estimate $|y_{j}'(t)|.$
Since $z_{j}(t)=m\phi^{m-1}(t+t_j)y_j'(t)$ solves the initial value problem
$z_{j}(s_{j}-t_{j})= w'(s_{j})/w(t_{j})\in[-1,0]$ for equation
$$
Dz'(t)-c\frac{1}{m\phi^{m-1}(t+t_j)}z(t)-p(t+t_j)y_j(t)-k(t+t_j)y_j(t-h)=0,
$$
we obtain that
\begin{align} \nonumber
z_{j}(t)=&e^{\frac{1}{D}\int_{s_j-t_j}^t\frac{c}{m\phi^{m-1}(\tau+t_j)}d\tau}
z_{j}(s_j-t_j)
\\  \label{eq-ode}
&+\frac{1}{D}
\int_{s_j-t_j}^t(p(t+t_j)y_j(s)+k(s+t_j)y_j(s-h))e^{\frac{1}{D}\int_s^t
\frac{1}{m\phi^{m-1}(\tau+t_j)}d\tau}ds.
\end{align}
In consequence,
\begin{equation}\label{eq-ybound}
|y_j'(t)| \le C+C(2|g'(\kappa)|+1) d,
\quad t \in[-2 h, s_{j}-t_{j}], \qquad j \in \mathbb{N},
\end{equation}
from which the uniform boundedness of the sequence $\{y_{j}'(t)\}$
on each compact interval $[-2 h, \xi], \xi>-2 h,$ follows.
Together with $0<y_{j}(t) \le \delta, t \ge-3 h,$ inequality \eqref{eq-ybound}
implies the pre-compactness of the set $\{y_{j}(t), t \ge-2 h, j \in \mathbb{N}\},$ in the compact open topology of $C([-2 h,+\infty), \mathbb{R}).$
Therefore, by the Arzela-Ascoli theorem combined with the diagonal method,
we can indicate a subsequence $y_{j_k}(t)$ converging to
a continuous function $y(t)$, $t \in[-2 h,+\infty).$
This convergence is uniform on every bounded subset of $[-2 h,+\infty).$
Additionally we may assume that
$\lim_{k \rightarrow \infty} y_{j_{k}}'(0)=y_{0}'$ exists.

Next, putting $s_{j}-t_{j}=0$ in \eqref{eq-ode}, we find that
\begin{align*}
z_j(t)=&m\phi^{m-1}(t+t_j)y_{j}'(t)\\
=&e^{\frac{1}{D}\int_{0}^t\frac{c}{m\phi^{m-1}(\tau+t_j)}d\tau}
z_j(0)
\\
&+\frac{1}{D}
\int_{0}^t(p(t+t_j)y_j(s)+k(s+t_j)y_j(s-h))e^{\frac{1}{D}\int_s^t
\frac{1}{m\phi^{m-1}(\tau+t_j)}d\tau}ds,~ t\ge -h.
\end{align*}
Integrating this relation between $0$ and $t$ and then taking the limit
as $j \rightarrow \infty$ in the obtained expression, we obtain that
\begin{align*}
y(t)=&1+\frac{Dm\kappa^{m-1}}{c}\left(e^{\frac{ct}{Dm\kappa^{m-1}}}-1\right)y_0'
\\
&+\int_0^t\frac{1}{Dm\kappa^{m-1}}\int_0^\sigma
(d'(\kappa)y(s)-b'(\kappa)y(s-h))e^{\frac{c(t-s)}{Dm\kappa^{m-1}}}
dsd\sigma, ~ t \ge-h.
\end{align*}

Therefore, $y(t)$ satisfies
\begin{equation}\label{eq-lineary}
Dm\kappa^{m-1}y^{\prime \prime}(t)-cy'(t)
-d'(\kappa)y(t)+b'(\kappa) y(t-h)=0, \quad t \ge-h.
\end{equation}
Additionally, $y(0)=1, y'(0)=y_{0}' \in[-1,0]$ and
$0 \le y(t) \le\delta, t \ge-2 h.$
Clearly, $y \in C^{2}(\mathbb{R}_{+})$
and we claim that $y(t)>0$ for all $t \ge0.$
Observe here that $y(t), t \ge-2 h,$ is non-increasing, and
therefore $y(0)=1, y(s)=0$ imply $s>0.$
Let us suppose, for a moment, that $y(s)=0$ and
$y(\tau)>0, \tau \in[-h, s).$
Then $y'(s)=0, y(s-h)>0,$ so that \eqref{eq-lineary} implies $y^{\prime \prime}(s)>0 .$ Thus $y(t)>0=y(s)$ for all $t>s$ close to $s$ which is not possible because $y$ is non-increasing on $[-2 h,+\infty)$.

We have proved that \eqref{eq-lineary} has a bounded positive solution on $\mathbb{R}_{+}$.
As it was established in \cite{Hupke} Lemma 3.1.1, this solution does not
decay superexponentially.
From Proposition 7.2 in \cite{Mallet} (see also Proposition 2.2 in \cite{Hale}),
we conclude that there are $b \le 0, \delta>0$
and a nontrivial eigensolution $v(t)$ of \eqref{eq-lineary}
on the generalized eigenspace associated with
the (nonempty) set $\Lambda$ of eigenvalues with $\Re \lambda=b,$ such that $y(t)=v(t)+O(\exp ((b-\delta) t))$, $t \rightarrow+\infty .$

On the other hand, since $c>c^{*},$ we know from Lemma \ref{le-characteristick}
that there are no real negative eigenvalues of \eqref{eq-lineary}
hence $\Im \lambda \neq 0$ for all $\lambda \in \Lambda .$
From Lemma 2.3 in \cite{Hale}, we find
that $y(t)$ is oscillatory, a contradiction.
$\hfill\Box$

{\it\bfseries Proof of Theorem \ref{th-osc}.}
This theorem follows from Lemma \ref{le-characteristick}
and Lemma \ref{le-notmonotone}.

Therefore, if $b'(\kappa)<0$ and the birth rate function $b$ satisfies
the feedback condition \eqref{eq-feedback},
then for $c>c_\kappa$, the semi-wavefront $\phi(t)$ is
slowly oscillating around the positive steady state.
In the remaining part of this section,
we show that these oscillations are non-decaying for $c$
greater than some constant $c^*$.

Before going further,
it will be convenient to work with the scaled function $\varphi(s)=\phi(cs)$.
Then $\varphi$ is a positive solution of the delay differential equation
$$
D\sigma(\varphi^{m})''(t)-\varphi'(t)-d(\varphi(t))+b(\varphi(t-r))=0,
\quad t \in \mathbb{R},
$$
where $\sigma=c^{-2}$.
The characteristic equation around $\kappa$ is
\begin{equation} \label{eq-character-star}
\chi^*(\lambda)=D\sigma m\kappa^{m-1}\lambda^2-
\lambda-d'(\kappa)+b'(\kappa)e^{-\lambda r}.
\end{equation}
We recall the following definition and auxiliary lemma in \cite{TrofimchukJDE}
concerned with the non-decaying oscillation around $\kappa$.

\begin{definition}[\cite{TrofimchukJDE}] \label{de-oscwavespeedc}
Suppose that $b'(\kappa)\le0$.
Let $c^*=c^*(m,r,b'(\kappa),d'(\kappa))\in(0,+\infty]$
be the largest extended real number such that
$\chi^*(\lambda)$ does not have roots in the half-plane $\{\Re z>0\}$
other than a positive real root.
\end{definition}

\begin{lemma} \label{le-cstar}
The inequality
$c^*(m,r,b'(\kappa),d'(\kappa))\ge c_\kappa(m,r,b'(\kappa),d'(\kappa))$
holds for all cases.
If $b'(\kappa)\ge -d'(\kappa)$, then $c^*(m,r,b'(\kappa),d'(\kappa))=+\infty$
for large time delay $r$;
while if $b'(\kappa)<-d'(\kappa)$, then
$$
c^*(m,r,b'(\kappa),d'(\kappa))=
\frac{\mu^*(m,b'(\kappa),d'(\kappa))+o(1)}{r},
\quad r\to+\infty,
$$
where $\mu^*(m,b'(\kappa),d'(\kappa)):=
\pi\sqrt{\frac{Dm\kappa^{m-1}}{-b'(\kappa)-d'(\kappa)}}$.
\end{lemma}
{\it\bfseries Proof.}
According to Lemma \ref{le-characteristick} and Lemma 1.1 in \cite{Gomez},
for any $c\le c_\kappa$, any complex root $\lambda_j$ of \eqref{eq-character-k}
has negative real part such that $\Re \lambda_j<\lambda_2<0$.
It follows that $c^*\ge c_\kappa$ for all cases.
If $c^*<+\infty$ and $c>c^*$, then \eqref{eq-character-star} has
a complex root with non-negative real part, denoted by $\alpha+i\beta$
with $\alpha\ge0$ and $\beta>0$.
Then
$$
Dm\kappa^{m-1}\sigma(\alpha+i\beta)^2
-(\alpha+i\beta)-d'(\kappa)+b'(\kappa)e^{-r(\alpha+i\beta)}=0.
$$
That is,
\begin{equation} \label{eq-zcstar}
\begin{cases}
Dm\kappa^{m-1}\sigma(\alpha^2-\beta^2)-\alpha-d'(\kappa)
+b'(\kappa)e^{-r\alpha}\cos(r\beta)=0,\\
2Dm\kappa^{m-1}\sigma \alpha \beta-\beta-b'(\kappa)e^{-r\alpha}\sin(r\beta)=0.
\end{cases}
\end{equation}
For large time delay $r$, we assert that $\alpha=o(1)$ as $r\to+\infty$.
Otherwise,
$|b'(\kappa)e^{-\lambda r}|<|D\sigma m\kappa^{m-1}\lambda^2-
\lambda-d'(\kappa)|$ for $\lambda\in \partial K$ for large time delay since
the complex-valued function
$D\sigma m\kappa^{m-1}\lambda^2-\lambda-d'(\kappa)$
has at most one complex root within $K$,
where $K:=\{z;\Re z>\alpha/2\}$ in the complex plane.
According to the Rouche's theorem,
\eqref{eq-character-star} admits at most one complex root
(that is a positive real number),
which is a contradiction.
Now, we see that
$$b'(\kappa)e^{-r\alpha}\cos(r\beta)=
d'(\kappa)+Dm\kappa^{m-1}\sigma\beta^2+\alpha
-Dm\kappa^{m-1}\sigma\alpha^2>d'(\kappa)+Dm\kappa^{m-1}\sigma\beta^2,$$
which is impossible if $b'(\kappa)\in[-d'(\kappa),0)$.
For the case $b'(\kappa)<-d'(\kappa)$,
we let $c$ tend to $c^*$, then $\alpha+i\beta$ tends to a purely imaginary number $iy$,
and the following limiting equation of \eqref{eq-zcstar} has a nonnegative solution
\begin{equation} \label{eq-zr}
\begin{cases}
-Dm\kappa^{m-1}\sigma y^2-d'(\kappa)
+b'(\kappa)\cos(ry)=0,\\
-y-b'(\kappa)\sin(ry)=0.
\end{cases}
\end{equation}
We note that according to the definition, $c^*$ is smallest positive real number
such that \eqref{eq-character-star} has complex roots with non-negative real part
except for the unique positive real root.
That is, $\sigma=1/(c^*)^2$ is the largest positive real number
such that \eqref{eq-zr} has a solution.
Asymptotic analysis as $r\to+\infty$ shows that
$ry\to\pi$ and $Dm\kappa^{m-1}\pi^2/(c^*r)^2\to -b'(\kappa)-d'(\kappa)$.
The proof is completed.
$\hfill\Box$

\begin{lemma}[\cite{TrofimchukJDE}]\label{le-stripeigenvalue}
If $c^*>0$ as in Definition \ref{de-oscwavespeedc},
then $\chi^*(\lambda)$ does not have
any zero in the strip $S_{00}:=(-\infty,0]\times
[-2\pi/r,2\pi/r]$ for every $c> c^*$.
\end{lemma}

Finally, similar to the proof of Theorem 3 in \cite{TrofimchukJDE},
we present a sufficient condition for the existence of non-decaying
oscillating semi-wavefronts.

\begin{lemma} \label{le-nondecayosci}
Assume that $b'(\kappa)<0$ and the birth rate function $b$ satisfies
the feedback condition \eqref{eq-feedback}.
If $c>c^*$, then the semi-wavefront $\phi(t)$
does not converge to $\kappa$ as $t\to +\infty$.
\end{lemma}
{\it\bfseries Proof.}
Using the similarly approach in \cite{TrofimchukJDE},
we can prove that the solution does not converge to $\kappa$, which implies that
the oscillation is non-decaying.
By contradiction, we assume that $\phi(t)\to \kappa$ as $t\to+\infty$.
Then $v(t)=\phi(t)-\kappa$ with $v(+\infty)=0$, satisfies
\begin{equation} \label{eq-initialv}
Dm\sigma(\phi(t)^{m-1}v'(t))'-v'(t)
-d_1(v(t))+b_1(v(t-r))=0,\quad t\in\mathbb R,
\end{equation}
where $b_1(s):=b(s+\kappa)-b(\kappa)$, $b_1(0)=0$, $b'(0)=b'(\kappa)$,
satisfies the feedback condition with respect to $0$,
and $d_1(s):=d(s+\kappa)-d(\kappa)$, $d_1(0)=0$, $d_1'(0)=d'(\kappa)$.

Since $v(+\infty)=0$, there exists a sequence
$t_{n} \rightarrow+\infty$ with the property such that
$\left|v(t_{n})\right|=\max _{s \ge t_{n}}|v(s)|$.
We can assume that $v$ attains its local extremum at
$t_n$ so that $v'(t_n)=0$, $v''(t)v(t_n)\le 0$.
These relations and \eqref{eq-initialv} imply that
$v(t_n)v(t_n-r)<0$ and therefore
$\operatorname{sc}(\overline v_{t_n})$ must be an odd integer.
Since $\operatorname{sc}(\overline v_{t_n})\le2$,
$\operatorname{sc}(\overline v_{t_n})=1$.
There are a unique $z_n\in(t_n-r,t_n)$
and a finite set $F_n$ such that $v(s)<0$ for
$s\in [t_n-r,z_n)\backslash F_n$ and $v(s)\ge 0$
for $s\in[z_n,t_n]$.
We can assume that
$|v(t_n)|=\max\{|v(s)|:s\in[z_n,t_n]\}$,
and that $\{r_n\}$, $r_n:=t_n-z_n\in(0,r)$,
is monotonically converging to $r^*\in[0,r]$.
Set $y_n(t)=v(t+z_n)/v(t_n)$, $t\in\mathbb R$, then $y_n(t)$ satisfies
$$
Dm\sigma(\phi(t)^{m-1}y'(t))'
-y'(t)-q_n(t)y(t)+p_n(t-h)y(t-h)=0,\quad t\in\mathbb R,
$$
where
$$
p_{n}(t)=
\begin{cases}
b_1\left(v\left(t+z_{n}\right)\right)/v\left(t+z_{n}\right),
&{\text { if } v\left(t+z_{n}\right) \neq 0}, \\
b'(\kappa),
&{\text { if } v\left(t+z_{n}\right)=0},
\end{cases}
$$
and
$$
q_{n}(t)=
\begin{cases}
d_1\left(v\left(t+z_{n}\right)\right)/v\left(t+z_{n}\right),
&{\text { if } v\left(t+z_{n}\right) \neq 0}, \\
d'(\kappa),
&{\text { if } v\left(t+z_{n}\right)=0}.
\end{cases}
$$
Since $y_{n}(0)=0$ and $\left|y_{n}(t)\right| \le1$, $t \ge0,$
and that $\lim_{n \rightarrow \infty} p_{n}(t)=b'(\kappa)$,
$\lim_{n \rightarrow \infty} q_{n}(t)=d'(\kappa)$,
$\lim_{n\rightarrow \infty}\phi(t)=\kappa$ uniformly in $t \in \mathbb{R}_{+}.$
From \eqref{eq-ode}, we get $|y_n(t)|$ is uniformly bounded in $C^1([-2r,\infty))$.
Hence, using the similar arguments in Lemma \ref{le-notmonotone},
there exists a sub-sequence $y_{n_j}$ converging to
$y^*(t)$, which is the solution of the linear equation
\begin{equation} \label{eq-lineary1}
D m\sigma\kappa^{m-1}y^{\prime \prime}(t)-y'(t)
-d'(\kappa)y(t)+b'(\kappa) y(t-h)=0, \quad t \ge2r.
\end{equation}
From Proposition 7.2 in \cite{Mallet},
for every sufficiently large $|\nu|$, $\nu<0$,
it holds that
$$
y^*(t)=Y_0(t)+O(\exp(\nu t)), \quad t \rightarrow+\infty,
$$
where $Y_0(t)$ is a nonempty finite sum of
eigensolutions of the linear equation
\eqref{eq-lineary1} associated to the eigenvalues in
$\{\lambda \in \mathbb{C}:\Re(\lambda) \in(-\nu, 0]\}$.
Thus, there exist $A>0, \beta>0, \alpha \ge 0, \zeta \in \mathbb{R},$ such that
$y^*(t)=(A \cos (\beta t+\zeta)+o(1)) e^{-\alpha t}, \quad t \ge2r$.
From Lemma \ref{le-stripeigenvalue} on the location of eigenvalues,
we have $\beta>2\pi/r$.
Since $y_{n_j}$ converges to $y^*$ as $j\to\infty$, this ensures that
$y_{n_j}$ changes its sign at least three times for sufficient large $j$.
It contradicts to $\operatorname{sc}(\overline v_{t_n})=1$
and completes the proof.
$\hfill\Box$

{\it\bfseries Proof of Theorem \ref{th-nondecay}.}
This theorem follows from Lemma \ref{le-cstar}, Lemma \ref{le-stripeigenvalue}
and Lemma \ref{le-nondecayosci}.

\section*{Acknowledgement}

This work was done when T.Y. Xu and S.M. Ji visited McGill University
supported by CSC programs.
They would like to express their sincere thanks for the hospitality
of McGill University and CSC.
The research of S. Ji was
supported by NSFC Grant No.~11701184 and CSC No. 201906155021
and the Fundamental Research Funds for the Central Universities of SCUT.
The research of M. Mei was supported in part
by NSERC Grant RGPIN 354724-16, and FRQNT Grant No. 2019-CO-256440.
The research of J. Yin was supported in
part by NSFC Grant No. 11771156 and NSF of Guangzhou Grant No. 201804010391.

\end{document}